\documentclass{article}

\usepackage[latin1]{inputenc}
\usepackage[square,numbers]{natbib}
\usepackage[english]{babel}

\makeindex

\everymath{\displaystyle} 
 
\usepackage{calc}
 
\usepackage{lmodern} 
\usepackage{theorem,latexsym,amsmath,amsfonts,amssymb,mathrsfs}
\usepackage[T1]{fontenc}
\usepackage{fancybox}
\usepackage{fancyhdr} 
\usepackage[english]{babel} 
\usepackage{multicol,multirow}
\usepackage{float}
\usepackage{array}
\usepackage[normalem]{ulem}
\usepackage{color}
\usepackage{tablists}
\usepackage{subfloat}
\usepackage{ragged2e}
\usepackage{graphicx}
\usepackage{geometry}
\usepackage{indentfirst}
\usepackage{latexsym}

\usepackage{sectsty} 
\usepackage{titlesec} 
\usepackage{enumerate} 
\usepackage{cancel} 
\usepackage{subfloat} 
\usepackage{sectsty} 
\usepackage[colorlinks, citecolor=blue,linkcolor=black]{hyperref}

\titleformat*{\subsection}{\bfseries\rmfamily}
\titleformat*{\section}{\bfseries\rmfamily}

\theorembodyfont{\slshape}

\usepackage[para]{footmisc}

\titlelabel{\thetitle.~} 
\numberwithin{equation}{section}

\AtBeginDocument{}

\newtheorem{Lemma}{Lemma}[section]{\bfseries}{\itshape}
\newtheorem{Proposition}{Proposition}[section]{\bfseries}{\itshape}

\newtheorem{Remark}{Remark} 

{\bfseries}{\itshape}
 
\newtheorem{Theorem}{Theorem}

\newcommand{\R}{\mathbb R}

\newcommand{\N}{\mathbb N}

\newcommand{\E}{\mathbb E}

\newcommand{\indep}{\perp\!\!\!\perp}

\title{
\textbf{Stochastic approximation method for kernel sliced average variance estimation}}
\date { }
\author{Emmanuel De Dieu  Nkou\\
\textit{Laboratoire URMI, Universit\'e des Sciences et Techniques de Masuku, Franceville, Gabon }\\
\textit{e-mail: emmanueldedieunkou@gmail.com} \\
}

\begin{document}
\maketitle

\begin{abstract}
In this paper, we use the stochastic approximation method to estimate Sliced Average Variance Estimation (SAVE). This method is known for its efficiency in recursive estimation. Stochastic approximation is particularly effective for constructing recursive estimators and has been widely used in density estimation, regression, and semi-parametric models. We demonstrate that the resulting estimator is asymptotically normal and root $n$ consistent. Through simulations conducted in the laboratory and applied to real data, we show that it is faster than the kernel method previously proposed.
\end{abstract}

\medskip 


\medskip 


\bigskip

\section{Introduction}\label{section1}

We consider the pair of random variables $(Y, X)$, where $Y$ is the response random variable and $X=\left(X_1,\cdots,X_d\right)^T$ the $d$-dimensional predictor random vector. The objective is to identify the optimal relationship between $X$ and $Y$. In this context,  one replaces $X$ by a linear combination of its components $B^TX$, $B =(\beta_1, . . . , \beta_N)$ is $d\times N$ matrix. $N \in \N^{*}$ such that $N <d$. This last condition is crucial for reducing the dimension of the subspace spanned by $\beta_1^TX, \cdots , \beta_N^TX$, known as the Effective Dimension Reduction (EDR) space.

To ensure that substituting $X$ with $B^TX$ does not result in any loss of information, it is necessary for  $Y$ and $X$  to be independent given $B^TX$(see \cite{li1991, duan1991, zhu1996, nkiet2008} among many others).  This condition can be expressed as
\begin{equation}\label{eqintro1}
Y\indep X | B^TX,
\end{equation} 
where ''$\indep$'' stands for independence. $B^T X$ summarizes the information in the predictors relevant to predicting $Y$. It is clear that the $EDR$ space is not unique. For uniqueness, we are interested in a subspace with minimal dimensions. Under mild conditions, this minimal subspace is often uniquely defined in practice and coincides with the intersection of all subspaces that satisfy condition (\ref{eqintro1}) \cite{cook1994, cook1996}.The intersection referred to as the central dimension-reduction ($CDR$) space is denoted as $S_{Y | X}$. In this article, we assume the existence of $S_{Y | X}$. Several methods are available in the literature to estimate the $CDR$ space, including $SIR$ (Sliced Inverse Regression) (see \cite{li1991, duan1991}) and $SAVE$ (Sliced Average Variance Estimation) (see \cite{cook1991, cook2000, li2007}). In these cited articles, the estimation of $B$ involves dividing the support of $Y$ into slices. However, alternative methods based on kernel estimation \cite{zhu1996, zhu2007} and wavelet estimation \cite{nkou2022a} have proven to be more effective. In this paper, our focus is on estimating the $CDR$ space using the $SAVE$ method with a kernel approach. We propose a novel recursive approach which, to our knowledge, has not been explored before.

In recursive estimation, the estimator can be updated with each new additional observation. This approach offers several advantages: the estimators are straightforward to implement, computationally efficient, and do not require substantial data storage. In specific cases, recursive estimators also appear to be more efficient than non-recursive counterparts.

This form of estimation was introduced by \cite{wolverton1969} and has evolved significantly in recent decades. Many references discuss this topic, including \cite{devroye1979} and \cite{amiri2009}. Among the widely used methods for recursive estimation is the stochastic approximation algorithm. This recursive method was initially introduced by \cite{revesz1973, revesz1977} and later extended by \cite{mokkadem2009b} to estimate regression functions.

The recursive approach has also been applied to estimate the $CDR$ space, as demonstrated in \cite{nguyen2010}. Recently, \cite{nkou2022c} utilized this approximation method to estimate the same space using the $SIR$ method, known for its efficiency over traditional approaches \cite{zhu2007}. According to these authors, a recursive estimator constructed via stochastic approximation primarily depends on two parameters: bandwidth and stepsize. The optimal selection of these parameters can significantly enhance the estimator's performance.

According to the $SAVE$ method, the directions $\beta_1, . . . , \beta_N$ are characterized under certain conditions as eigenvectors of the $d \times d$ matrix $\Lambda=\E\big(I_d - cov(X|Y)\big)^2$. Therefore, the aforementioned estimation problem is based on estimating $\Lambda$. Building on the estimator of $\Lambda$ using the kernel method proposed by \cite{zhu2007}, we propose constructing one using the stochastic approximation method and investigate its asymptotic normality.

The remainder of the paper is structured as follows. In Section \ref{section2}, we present a recursive estimator constructed using recursive estimators for density and regression functions. Section \ref{section3} outlines our main results. Section \ref{section4} presents simulation results, where a real dataset is used to illustrate our approach. Section \ref{section5} contains the proofs of the main results. We conclude in Section \ref{section6} with proofs of additional important results that support the theorems.

\medskip

\section{Kernel sliced average variance estimation using the stochastic approximation method}\label{estimator}\label{section2}
Consider the pair of random variables $(Y, X)$, where $Y$ represents the response random variable and $X=\left(X_1,\cdots,X_d\right)^T$ denotes the 
$d-$dimensional predictor random vector. It is assumed, without loss of generality, that  $\E \left(X\right) = 0$. This implies, for instance, that there exists a random variable $Z$ in $\R^d$ such that its covariance matrix $\Sigma_Z = var(Z)$ is invertible, and   $X = \Sigma_Z^{-1/2}\left(Z-\E(Z)\right)$. The earlier equation yields the relationship $S_{Y |Z} = \Sigma_Z^{-1/2} S_{Y |X}$ (see \cite[Proposition 6.1]{cook1998}). The columns of the matrix $\Sigma_Z^{1/2}\eta$, where $\eta$ forms an appropriate basis of $S_{Y |Z}$, constitute a basis of the $CDR$ space of $Y$ on $X$. This explains why there is no loss of information when using such a transformation. Additionally, we assume $\E\left(X_j^2\right)<+ \infty$, $j=1,\cdots,d$.  Under certain conditions, the directions $\beta_1,\cdots, \beta_N$ are eigenvectors of the $SAVE$ matrix $\Lambda$ (see \cite{zhu2007}).

Considering  $A^2 = AA$ as a squared symmetric matrix $A$, then the $SAVE$ matrix is defined as
\begin{equation}\label{rs1}
\Lambda = \E\bigg(I_d - cov(X|Y)\bigg)^2=I_d-2\,\E\bigg(cov(X|Y)\bigg) + \E\bigg(cov(X|Y)\bigg)^2.
\end{equation}
$I_d$ denotes the identity matrix of order $d$ and $cov(X|Y)$ represents the conditional covariance matrix of  $X$ given $Y$. When using $SAVE$  to identify the  $CDR$ space  $E :=S_{Y |X}$, we need to assume the following two conditions:
\begin{equation}\label{cond1}
\E\left(X|P_{E}X\right)=P_{E}X \quad \mbox{ and }\quad cov\left(X|P_{E}X\right)= I_d - P_{E}. 
\end{equation}
$P_{E}$ stands the linear projection operator onto space $E$ with respect to the standard inner product in $\R^d$. If $X$ follows an elliptically contoured distribution, the first equality is automatically satisfied and second equality approximately holds. If $X$ follows a multivariate normal distribution, both equalities are guaranteed (see \cite{cook1991, cook1998}). In cases where these conditions are not satified, $X$ can often be transformed accordingly to induce them.

To define a kernel recursive estimator of $\Lambda$ using the stochastic approximation method, we introduce additional functions. For $j,\,k,\,\ell \in \{1,\cdots,d\}$, we consider the functions $r_j$ and $R_{k\ell}$ defined on $\R$ respectively by:
$$
r_j(y)=\E\left(X_j|Y=y\right)\quad \mbox{ and }\quad R_{k\ell}(y)=\E\left(X_k X_\ell|Y=y\right).
$$
In what follows, we will consider $\delta_{k\ell}$ which is such that $\delta_{k\ell}= 1$ if $k = \ell$ and $\delta_{k\ell} = 0$ otherwise.\\
Setting \,$\Lambda = \big(\lambda_{k\ell}\big)_{1\leqslant k,\ell \leqslant d}$,\, we have
\begin{multline*}
\lambda_{k\ell}=\delta_{k \ell}-2\,\E\bigg(R_{k\ell}(Y) - r_k(Y)r_\ell(Y)\bigg)\\ + \E\Bigg(\sum_{j=1}^d\bigg(R_{kj}(Y)R_{j,\ell}(Y)-R_{kj}(Y)r_j(Y)r_\ell(Y)\\
 - r_k(Y)r_j(Y)R_{j\ell}(Y) + r_k(Y)r_{\ell}(Y) r_j^2(Y)\bigg)\Bigg).
\end{multline*}
It is clear that from the above,  it is evident that to obtain an estimate of $\lambda_{k,\ell}$ using a  recursive method  based on a stochastic approximation algorithm, we require  estimators of $r_j$ and $R_{k\ell}$ obtained through the  same method.  The estimator for $r_j$ has already been proposed and can be found in   \cite{slaoui2016} and in \cite{slaoui2020}.

Let $f$ denote the density function of $Y$. We assume that for all $y\in\R$, $f(y) > 0$. We consider
$$
g_j(y)=r_j(y) f(y) = \int x f_{_{(X_j, Y)}} (x,y)dx
$$
and 
$$
 G_{k\ell}(y)=R_{k\ell}(y) f(y) = \int\int  x_k x_\ell f_{_{(X_k , X_\ell, Y)}} (x_k,x_\ell,y)dx_k dx_\ell.
$$
Simultaneously, we introduce the functions  $v_j$ and $v_{k\ell}$ which will be useful to us later.
$$
v_j(y) = \int x^2 f_{_{(X_j, Y)}} (x,y)dx \, \mbox{ , }\, v_{jk\ell}(y) = \int x_k^2 x_\ell^2 f_{_{(X_{j,k}, X_{j,\ell}, Y)}} (x_k, x_\ell,y)dx_k dx_\ell.
$$
The estimators of $r_j$ and $R_{k\ell}$ are thus obtained by those of $f$, $g_j$ and $G_{k\ell}$.

Let \,$X_i = \left(X_{i,1}\,,\cdots,\,X_{i,d}\right)$. Consider an i.i.d. sample $\left(Y_i, X_i\right)_{i=1,...n}$ from the pair $(Y, X)$. From this sample of size $n$,  the kernel estimator $\bar{G}_{k\ell,n}$   of $G_{k\ell}$ is given by:
$$
 \bar{G}_{k\ell,n}(y)=\frac{1}{n}\sum_{i=1}^n X_{i, k}X_{i,\ell}\frac{1}{h_n}K\left(\frac{y-Y_i}{h_n}\right), 
$$ 
$K$ is a kernel and $(h_n)$ denotes a bandwidth. The kernel estimators of $f$ and $r_j$ are well known. From these  estimators, we obtain the kernel estimator of $\lambda_{k\ell}$ as described in \cite{zhu2007}.

The Monro's scheme \cite{revesz1973, revesz1977} for constructing approximation algorithms to find  the zero $x^*$ of an unknown function $h : \R \rightarrow \R$ proceeds as follows : initially, $x_0$ is arbitrarily chosen, and then the sequence $(x_n)$ is recursively defined by
$$
x_{n+1}=x_{n} + \gamma_{n+1} W_{n+1}, 
$$
where $W_{n+1}$ is an observation of the function $h$ at  point $x_{n}$. The three estimators are  directly derived using this scheme, which is a method introduced by \cite*{mokkadem2009a} for constructing recursive kernel estimators of a probability density.

Suppose we add an additional $(n+1)$th pair $\left(Y_{n+1}, X_{n+1}\right)$ of random variables. From this new data, we seek to obtain the estimators $\widehat{f}_{n+1}$, $\widehat{g}_{j,n+1}$ and $\widehat{G}_{k\ell,n+1}$ using the stochastic approximation method, we have:
\begin{equation}\label{eqnfgga1}
\begin{split}
\widehat{f}_{n+1}(y)        &= \, (1 - \gamma_{n+1})\widehat{f}_{n}(y) +  \frac{\gamma_{n+1}}{h_{n+1}} \, K\left(\frac{y-Y_{n+1}}{h_{n+1}}\right),\\
\widehat{g}_{j,n+1}(y)      &= (1 - \gamma_{n+1})\widehat{g}_{j,n}(y) +  \frac{\gamma_{n+1}}{h_{n+1}}\, X_{n+1,j} \,K\left(\frac{y-Y_{n+1}}{h_{n+1}}\right),\\
\widehat{G}_{k\ell,n+1}(y) &= \, (1 - \gamma_{n+1})\widehat{G}_{k\ell,n}(y) +  \frac{\gamma_{n+1}}{h_{n+1}}\, X_{n+1,k} \, X_{n+1,\ell}\, K\left(\frac{y-Y_{n+1}}{h_{n+1}}\right).
\end{split}
\end{equation}
The stepsize $(\gamma_n)$ is a sequence of positive real numbers that converges to zero.
Using these constructions, directly, we have
\begin{multline*}
\widehat{\lambda}_{k\ell,n+1} = \frac{1}{n+1}\sum_{i=1}^{n+1}\Bigg( \delta_{k \ell}-2\,\bigg(\widehat{R}_{k\ell,n+1}(Y_i) - \widehat{r}_{k,n+1}(Y_i)\widehat{r}_{\ell,n+1}(Y_i)\bigg)\\ + \sum_{j=1}^d\bigg(\widehat{R}_{kj,n+1}(Y_i)\widehat{R}_{j\ell,n+1}(Y_i)-\widehat{R}_{kj,n+1}(Y_i)\widehat{r}_{j,n+1}(Y_i)\widehat{r}_{\ell,n+1}(Y_i) \\
- \widehat{r}_{k,n+1}(Y_i)\widehat{r}_{j,n+1}(Y_i)\widehat{R}_{j\ell,n+1}(Y_i) + \widehat{r}_{k,n+1}(Y_i)\widehat{r}_{\ell,n+1}(Y_i) \widehat{r}_{j,n+1}^2(Y_i)\bigg)\Bigg).
\end{multline*}
It is worth noticing that if we set $\widehat{f}_{0}(y) = \widehat{g}_{j,0}(y) = \widehat{G}_{k\ell,0}(y)=0$, then for $n\geq 1$,
\begin{equation}\label{eqnfgg1}
\begin{split}
\widehat{f}_{n}(y)        &= \, \pi_n \sum_{i=1}^n \pi_i^{-1} \gamma_i \frac{1}{h_i} K\left(\frac{y-Y_i}{h_i}\right),\\
\widehat{g}_{j,n}(y)      &= \, \pi_n \sum_{i=1}^n \pi_i^{-1} \gamma_i X_{i,j} \frac{1}{h_i}K\left(\frac{y-Y_i}{h_i}\right),\\
\widehat{G}_{k\ell,n}(y)  &= \, \pi_n \sum_{i=1}^n \pi_i^{-1} \gamma_i X_{i,k} X_{i,\ell} \frac{1}{h_i}K\left(\frac{y-Y_i}{h_i}\right),
\end{split}
\end{equation}
where the sequence $(\pi_n)$ is defined by:
\begin{equation}\label{eqnpi1}
 \pi_0=1\, \mbox{ and for  }\, n \geq 1, \,\, \pi_n = \prod_{i=1}^n (1-\gamma_{i}),
\end{equation}
we derive recursive estimators of $r_j$ and $R_{k\ell}$ by
\begin{equation}\label{eqnsr1}
\widehat{r}_{j,n}(y) = \frac{\widehat{g}_{j,n}(y)}{\widehat{f}_{n}(y)}\quad \mbox{ and }\quad \widehat{R}_{k\ell,n}(y)=\frac{\widehat{G}_{k\ell,n}(y)}{\widehat{f}_{n}(y)},
\end{equation}
we  obtain an estimator of $\lambda_{k\ell}$  whose formulation resembles  that in \cite{zhu2007}.
Nevertheless, with a specific set of operations, we can derive the following recursive form:
\begin{equation}\label{eqnasr2}
\widehat{\lambda}_{k\ell,n+1} = \frac{n}{n+1}\widehat{\lambda}_{k\ell,n} + \mathcal{W}_{k\ell,n+1}.
\end{equation}
The sequence $\mathcal{W}_{k\ell,n+1}$ involves a combination of addition and multiplication of 
$\widehat{f}_{n}(y)$, \, $\frac{\gamma_{n+1}}{h_{n+1}} \, K\left(\frac{y-Y_{n+1}}{h_{n+1}}\right)$,\, $\widehat{g}_{j,n}(y)$, $\frac{\gamma_{n+1}}{h_{n+1}}\, X_{n+1,j} \,K\left(\frac{y-Y_{n+1}}{h_{n+1}}\right)$,\, $\widehat{G}_{k\ell,n}(y)$,\, $\frac{\gamma_{n+1}}{h_{n+1}}\, X_{n+1,k} \, X_{n+1,\ell}\, K\left(\frac{y-Y_{n+1}}{h_{n+1}}\right)$.
\begin{Remark}
In Equations (\ref{eqnfgga1}) and (\ref{eqnasr2}),  this corresponds to adding one additional datum. These equalities can be extended to the case where $p$ data points are added. We refer to Remark 2.2 of \cite{nkou2022c}. However, generalizing equation (\ref{eqnasr2}) to accommodate $p$ data points would be quite cumbersome. This approach was not pursued.
\end{Remark}
To avoid division by zero, we will consider the following transformations, as proposed in \cite{zhu1996}.
$$
f_{b_n}(y)=\max\big(f(y),b_n\big)\,\,\textrm{ , }\,\, \, \widehat{f}_{b_n}(y)=\max\left(\widehat{f}_n(y),b_n\right),
$$
$$
R_{k\ell,b_n}(y)=\frac{G_{k\ell}(y)}{f_{b_n}(y)}\,\,\textrm{ and }\,\,\, r_{j,b_n}(y) = \frac{g_j(y)}{f_{b_n}(y)}, \quad j=1,\cdots,d.  
$$

The sequence $\left(b_n\right)_{n\in\mathbb{N}^\ast}$ consists of positive real numbers that satisfy the property:  $\lim_{n\rightarrow +\infty}\, b_n =0$.

Consequently, we define
\begin{multline}\label{lambda_bn1}
\lambda_{k\ell,b_n}=\delta_{k \ell}-2\,\E\bigg(R_{k\ell,b_n}(Y) - r_{k,b_n}(Y)r_{\ell,b_n}(Y)\bigg)\\ + \E\Bigg(\sum_{j=1}^d\bigg(R_{kj,b_n}(Y)R_{j\ell,b_n}(Y)-R_{kj,b_n}(Y)r_{j,b_n}(Y)r_{\ell,b_n}(Y)\\
 - r_{k,b_n}(Y)r_{j,b_n}(Y)R_{j\ell,b_n}(Y) + r_{k,b_n}(Y)r_{\ell,b_n}(Y) r_{j,b_n}^2(Y)\bigg)\Bigg),
\end{multline}
and we have the estimators :
$$
\widehat{r}_{j,b_n}(y) = \frac{\widehat{g}_{j,n}(y)}{\widehat{f}_{b_n}(y)}\quad \mbox{ and }\quad \widehat{R}_{k\ell,b_n}(y)=\frac{\widehat{G}_{k\ell,n}(y)}{\widehat{f}_{b_n}(y)}.
$$
A  recursive kernel  estimator of $\lambda_{k\ell}$ can therefore be formulated as follows:
\begin{multline}\label{eqnest1}
\widehat{\lambda}_{k\ell,b_n} = \frac{1}{n}\sum_{i=1}^n\Bigg( \delta_{k \ell}-2\,\bigg(\widehat{R}_{k\ell,b_n}(Y_i) - \widehat{r}_{k,b_n}(Y_i)\widehat{r}_{\ell,b_n}(Y_i)\bigg)\\ + \sum_{j=1}^d\bigg(\widehat{R}_{kj,b_n}(Y_i)\widehat{R}_{j\ell,b_n}(Y_i)-\widehat{R}_{kj,b_n}(Y_i)\widehat{r}_{j,b_n}(Y_i)\widehat{r}_{\ell,b_n}(Y_i) \\
- \widehat{r}_{k,b_n}(Y_i)\widehat{r}_{j,b_n}(Y_i)\widehat{R}_{j\ell,b_n}(Y_i) + \widehat{r}_{k,b_n}(Y_i)\widehat{r}_{\ell,b_n}(Y_i) \widehat{r}_{j,b_n}^2(Y_i)\bigg)\Bigg).
\end{multline}

Note that with the introduction of $b_n$, we obtain for all $y\in \R$:
\begin{equation}\label{eqnbn1}
|\widehat{f}_{b_n}(y) - f_{b_n}(y)|\leq |\widehat{f}_{n}(y) - f(y)|\quad \mbox{ , } \quad |R_{k\ell,b_n}(y)| \leq |R_{k\ell}(y)|\quad \mbox{ , } \quad |r_{j,b_n}(y)| \leq |r_{j}(y)|.
\end{equation}

\medskip

\section{Assumptions and main results}\label{section3}
In this section, we present our main result on the asymptotic distribution of $\sqrt{n}\left(\widehat{\Lambda}_n  - \Lambda\right)$. 
Recall that we have an independent, identically distributed sample $(Y_i, X_i)_{i=1,...,n}$ of the pair of random variables $(Y, X)$ satisfying the model (\ref{eqintro1}). To establish this result, we will require the following assumptions:

\begin{description}
\item[A1.] For any sample $X_1,..., X_n$ of independent, identical random variables with the same distribution of $X$, there exists a sequence $M_n$ of strictly positive numbers such that $M_n\sim \sqrt{\log n}$ and  satisfies 
$$
\sup_{1\leq i \leq  n}  \left\|X_i\right\|_d \leq M_n,
$$
where $\left\|\cdot \right\|_d$ denotes the usual euclidean norm in $\R^d$.
\end{description}
Assumption \textbf{A1.} holds true even when $X$ is bounded. This condition is particularly satisfied when $X$ follows a multivariate distribution.

\begin{description}
\item[A2.]  The function $H$ represents  one of the following functions: the density $f$ of $Y$, the functions $g_j = f R_j$ and the functions $G_{k\ell}$. $H$ satisfies the following conditions:
$H \in \mathcal{L}^2(\R)$is  three times differentiable and its third  derivatives satisfy the Lipschitz property: there exists $c>0$ such that
$$
\left|H^{(3)}(y+u) - H^{(3)}(y)\right|\leq  c|u|,
$$
for $1\leq j, k,\ell \leq d$. The constant $c$ may vary in different contexts. 
\end{description}

\begin{description}
\item[A3.] The function $H$ represents  one of the following functions : $v_j$,  $r_j$, $R_{k\ell}$, $r_k r_\ell$, $R_{kj}R_{j\ell}$, $R_{k\ell}r_j$, $R_{k\ell}r_\ell r_j$, $r_k r_\ell r_j$ for $1\leq k,\ell, j \leq d$.

$H$ is such that $\int \left|H(y)\right|dy < + \infty$ and there exists $c>0$ such that \,$\left|H(y+u) - H(y)\right|\leq  c|u|$.
\end{description}

The conditions on $f$ and $g_j$ in Assumption \textbf{A.2} have been used in several papers addressing similar topics, specifically in \cite{zhu1996} and in \cite{zhu2007}. The latter utilized  Assumption \textbf{A.3}.
\begin{description}
\item[A4.] The kernel $K : \R\mapsto \R$ satisfies:
\begin{enumerate}[(1).]
	\item $K(\cdot)$ is symmetric about 0 , for all $u\in\R$, $K(u)>0$ and $\int_{\R} K(u)du = 1$;
	\item $D=\sup_{u\in \R}K(u)<+\infty$; 
	\item $\int_{\R} u^s K(u) du =0$ for $s\in \left\{1,2,3\right\}$ and $\int_{\R} u^4 K(u) du < +\infty$;
	\item $\int u K^2(u) du = 0$,  $\int  K^2(u) du< +\infty$ and $\int u^s K^2(u) du < +\infty$  for $s\in \left\{2,8\right\}$ .
\end{enumerate}
\end{description}
Assumption \textbf{A4} encompasses the standard conditions typically satisfied by most kernels.
\begin{description}
\item[A5.] The bandwidth $(h_n)$,  the sequence $(b_n)$ and  the stepsize $(\gamma_n)$ verify the following conditions:
\begin{enumerate}[(1).]
	\item $(h_n) \in \mathcal{GS}(-c_1)$ and $(b_n) \in \mathcal{GS}(-c_2)$, where  $c_1$ and $c_2$ are numbers satisfying: $1/5< c_1 <1/4 - 2c_2$, $1/50<c_2<1/25$;
	\item $(\gamma_n) \in \mathcal{GS}(-1)$.
\end{enumerate}
\end{description}
\begin{Remark}
Assumption \textbf{A5.} combines the assumptions on  $(h_n)$,  $(b_n)$  and  $(\gamma_n)$ used by \cite{slaoui2015, mokkadem2009a, zhu1996}. Specifically, $(\gamma_n)$ is chosen to minimize the variance of $\widehat{f}_{n}$. Regarding  $(b_n)$, we assume   $b_n=\min(\varepsilon,n^{-c_2})$, where $\varepsilon$ is a fixed strictly positive number chosen to be sufficiently small. This choice enhances the accuracy of the estimation of  $f$; see  further details in  Remark 3.1 in \cite{nkou2019}.  It is important to note that despite this consideration $(b_n) \in \mathcal{GS}(-c_2)$. Conditions on $(\gamma_n)$ are standard in the framework of stochastic approximation algorithms, implying in particular that the limit of $(n\gamma_n)^{-1}$  is finite. In the following, we will use the notation $\varepsilon_0$ which is such that:
\begin{equation}\label{eqga1}
\varepsilon_0 = \lim_{n\rightarrow + \infty}(n\gamma_n)^{-1}. 
\end{equation}
\end{Remark}

Note that $(v_n) \in \mathcal{GS}(v^*)$ if $\lim_{n\rightarrow + \infty} n\left(1-\frac{v_{n-1}}{v_n}\right) = v^*$. For more details on the sets $\mathcal{GS}$, refer to \cite{galambos1973} in a more general context and to \cite{mokkadem2007} in the context of stochastic approximation algorithms.
\begin{description}
\item[A6.] $\sqrt{n}\E\left[R_{k\ell}^2(Y)\mathbf{1}_{\{f(Y)<b_n\}}\right] = 0$ and $\sqrt{n}\E\left[r_{j}(Y)r_{k}(Y)r_{\ell}(Y)\mathbf{1}_{\{f(Y)<b_n\}}\right] = 0$ as $n\rightarrow + \infty$ for $1\leq j, k,\, \ell \leq d$, where $\mathbf{1}_{(\cdot)}$ is the indicator function and $(b_n)$ satisfies Assumption \textbf{A5.}.
\end{description}
Assumption  \textbf{A6.} is similar to Assumption 6 introduced by \cite{zhu1996}.
\medskip

We define the following random variables for $1 \leq j,\, k,\, \ell \leq d$:
\begin{eqnarray*}
\mathcal{A}^{(1)}_{k\ell}(X_i,Y_i)&=&\big(X_{i,k}X_{i,\ell}- R_{k\ell}(Y_i)\big)\\
\mathcal{A}^{(2)}_{k\ell}(X_i,Y_i)&=& \big(X_{i,k}r_{\ell}(Y_i) + X_{i,\ell}r_{k}(Y_i) - 2r_{k}(Y_i)r_{\ell}(Y_i)\big)\\
\mathcal{B}^{(1)}_{jk\ell}(X_i,Y_i)&=&\Big(X_{i,\ell}X_{i,j}R_{k j}(Y_i) + X_{i,k}X_{i,j}R_{j\ell}(Y_i)- 2R_{j\ell}(Y_i)R_{kj}(Y_i)\Big)\\
 \mathcal{B}^{(2)}_{(j\ell) k}(X_i,Y_i)&=&\bigg(X_{i,\ell}X_{i,j}r_{j}(Y_i)r_{k}(Y_i) + X_{i,j}R_{j \ell}(Y_i)r_{k}(Y_i)\\
 & & \hspace{5cm} +\, X_{k,i}R_{\ell j}(Y_i)r_{j}(Y_i) -3r_{k}(Y_i)R_{j\ell}(Y_i)r_{j}(Y_i)\bigg)\\
\mathcal{B}^{(3)}_{(j\ell) k}(X_i,Y_i)&=&\bigg(X_{i,\ell}X_{i,j}r_{j}(Y_i)r_{k}(Y_i) + X_{i,j}R_{\ell j}(Y_i)r_{k}(Y_i)\\
 & & \hspace{5cm} +\, X_{k,i}R_{j\ell}(Y_i)r_{j}(Y_i) -3r_{k}(Y_i)R_{j\ell}(Y_i)r_{j}(Y_i)\bigg)\\
\mathcal{B}^{(4)}_{j k\ell}(X_i,Y_i)&=&\bigg(X_{i,k}r_{j}^2(Y_i)r_{\ell}(Y_i) + X_{i,\ell}r_{j}^2(Y_i)r_{k}(Y_i)+2X_{i,j}r_{j}(Y_i)r_{k}(Y_i)r_{\ell}(Y_i)\\
  & & \hspace{5cm} - \, 4r_{j}^2(Y_i)r_{k}(Y_i)r_{\ell}(Y_i)\bigg),
\end{eqnarray*}
and
\begin{multline*}
\Gamma_{k\ell}(X_i,Y_i)=\delta_{k,\ell}-2\left(\mathcal{A}_{k\ell}^{(1)}(X_i,Y_i)-\mathcal{A}_{k\ell}^{(2)}(X_i,Y_i) \right)\\
+\,\sum_{j=1}^d\bigg(\mathcal{B}_{jk\ell}^{(1)}(X_i,Y_i) - \mathcal{B}_{(j\ell) k}^{(2)}(X_i,Y_i)- \mathcal{B}_{(kj) \ell}^{(2)}(X_i,Y_i) \\
- \mathcal{B}_{(j\ell) k}^{(3)}(X_i,Y_i)- \mathcal{B}_{(k j) \ell}^{(3)}(X_i,Y_i) + \mathcal{B}_{jk\ell}^{(4)}(X_i,Y_i)\bigg).
\end{multline*}
Considering the random variables below, we define the $k\ell$th element of the matrix $\Gamma_n$ as follows:
$$
\Gamma_{k\ell,n}=\pi_n \sum_{i=1}^n  \pi_i^{-1}\gamma_i\Gamma_{k\ell}(X_i,Y_i).
$$
Likewise, we define
\begin{eqnarray*}
T_{k\ell}(X_i,Y_i)&=& \delta_{k\ell} -2 \big(X_{i,k}X_{i,\ell}+X_{i,k}r_{\ell}(Y_i)+X_{i,\ell}r_{k}(Y_i)- r_{k}(Y_i)r_{\ell}(Y_i)\big)\\
                  &&+ \sum_{j=1}^d\big(X_{i,k}X_{i,j}R_{j\ell}(Y_i) + X_{i,\ell}X_{i,j}R_{jk}(Y_i) - R_{kj}(Y_i)R_{j\ell}(Y_i)\big)\\
									&&+\sum_{j=1}^d\big(X_{i,k}X_{i,j}r_{j}(Y_i)r_{\ell}(Y_i) + X_{i,j}r_{\ell}(Y_i)R_{kj}(Y_i) + X_{i,\ell}r_{j}(Y_i)R_{kj}(Y_i)\\
									&& \hspace{8cm} -\, 2r_{\ell}(Y_i)r_{j}(Y_i)R_{kj}(Y_i) \big)\\
									&&+\sum_{j=1}^d\big(X_{i,j}X_{i,\ell} r_{\ell}(Y_i) + X_{i,j}r_{\ell}(Y_i)R_{kj}(Y_i) + X_{i,\ell}r_{j}(Y_i)R_{kj}(Y_i)\\
									&& \hspace{8cm}-\, 2r_{j}(Y_i)r_{\ell}(Y_i)R_{kj}(Y_i) \big)\\
									&&+\sum_{j=1}^d\big(X_{i,k}r_{j}^2(Y_i)r_{\ell}(Y_i) + X_{i,\ell}r_{j}^2(Y_i)r_{j}(Y_i)+2X_{i,j}r_{j}(Y_i)r_{k}(Y_i)r_{\ell}(Y_i)\\
									&& \hspace{8cm}-\, 3r_{j}^2(Y_i)r_{k}(Y_i)r_{\ell}(Y_i)\big),
\end{eqnarray*}
and  the $k\ell$th element of the matrix $T_n$ as: \, $T_{k\ell,n}=\frac{1}{n}\sum_{i=1}^n T_{k\ell}(X_i,Y_i)$.

\medskip

We now present the main results.

\begin{Theorem}\label{theo1}
Under the  Assumptions \textbf{A1.} - \textbf{A5.}, as $n\rightarrow +\infty$, we have
\[
\sqrt{n}\left(\widehat{\Lambda}_n-\Lambda\right)\stackrel{\mathscr{D}}{\longrightarrow}\Phi,
\]
where $\stackrel{\mathscr{D}}{\longrightarrow}$ denotes convergence in distribution and 
 $\Phi$  is a random variable having a normal distribution,  in the space $\mathscr{M}_d(\R)$ of $d\times d$ matrix, such that, for any $A=\left(a_{k\ell}\right)\in \mathscr{M}_d(\mathbb{R})$, $A\neq 0$,  one has   $\textrm{tr}\Big(A^T \Phi\Big) \leadsto \mathcal{N}(0,\sigma^2_{A})$  is given by:  
\begin{equation}\label{siga}
\sigma^2_{A}=Var\bigg(\sum_{k=1}^d \sum_{\ell=1}^d a_{k\ell}\Big(\Gamma_{k\ell}(X,Y) + T_{k\ell}(X,Y)\Big)\bigg).
\end{equation}
\end{Theorem}

\medskip

A consequence of Theorem \ref{theo1} is  the asymptotic normality of the eigenvector estimators. For $(j,k)\in\{1,\cdots,d\}^2$, let  $\beta_j=\left(\beta_{j1},\cdots,\beta_{jd}\right)^T$ denote an eigenvector of  $\Lambda$ and $\widehat{\beta}_{j,n}=\left(\beta_{j1,n},\cdots,\beta_{jd,n}\right)^T$ denote its estimator, which is an eigenvectors of $\widehat{\Lambda}_n$. Let $\lambda_j$, a eigenvalue of $\beta_j$. We assume $\lambda_1> \cdots \lambda_d>0$. This assumption has been previously used in \cite{zhu1996, zhu2007}. We consider the random variable
$$
\mathcal{V}_{jk}=\left(\sum_{\stackrel{r=1}{r\neq j}}^d\frac{\beta_{rk}}{\lambda_j-\lambda_r}\right)\,\sum_{p=1}^d\sum_{q=1}^d\beta_{jp}\beta_{jq}\Big(\Gamma_{pq}(X,Y) + T_{pq}(X,Y)\Big).
$$
Thus, we have:

\begin{Theorem}\label{theo2}
Under the  Assumptions \textbf{A1.} - \textbf{A6.},  we have
$$
\sqrt{n}\left(\widehat{\beta}_{j,n}-\beta_j\right)\stackrel{\mathscr{D}}{\longrightarrow}\mathcal{N}(0,\Sigma_j)
\, , \mbox{ as }\, n\rightarrow +\infty,
$$
where $\Sigma_j$ is the $d\times d$ covariance matrix of the random vector $ \mathcal{V}_{j}=\left(\mathcal{V}_{j1},\cdots,\mathcal{V}_{jd}\right)^T$.
\end{Theorem}
 The proof of Theorem \ref{theo1} is presented in Section \ref{section5_2}, while that of Theorem \ref{theo2} follows a similar approach to  Theorem 2 in \cite{nkou2022a}; therefore, we omit the details in this paper. 

\medskip

\section{Simulation study and application to real data}\label{simulations1}\label{section4}
In this section, we evaluate the performance of the proposed recursive estimator ($Save-R$) on simulated laboratory samples and real data, comparing it with the non-recursive estimator proposed by \cite{zhu2007}($Save-NR$). This study is the main focus of the paper. $SAVE$ is more efficient than $SIR$ \cite{zhu2007}, it involves additional computational steps.

 Therefore, we do not compare the recursive estimator of $SIR$ \cite{nkou2022c} with the recursive estimator of SAVE presented in this paper. 

We assess the convergence speed of the estimators using the correlation coefficient $R^2$ below proposed by \cite{li1991}, which measures the squared cosine of the angle between two vectors. Specifically, it quantifies the fit between the estimated $CDR$ space $ \left(\widehat{\beta}_j\right)$ and the true $CDR$ space $S_{Y |X}$ $(\beta_j)$.

A coefficient closer to one indicates a better estimation, meaning that the estimated. That is $CDR$ space aligns well with the vectors $\left(\widehat{\beta}_1,\cdots, \widehat{\beta}_N\right)$  associated with the $N$ largest eigenvalues of $\widehat{\Lambda}_n$. This criterion is $R^2\left(\widehat{\beta}_j\right)$ defined by:
$$
R^2\left(\widehat{\beta}_j\right)=\frac{\left(\widehat{\beta}^{T}_j   \beta_j\right)^2}{\widehat{\beta}^{T}_j \widehat{\beta}_j\cdot \beta^T_j  \beta_j}.
$$
See \cite{li1991} for more details on $R^2\left(\widehat{\beta}_j\right)$. In terms of computational efficiency, we compare the time taken by each estimator to complete the task. The results were obtained using identical samples for a fair comparison.The simulations were conducted using R software. The corresponding codes are available from the author upon request. For the construction both estimators, we employed Epanechnikov kernel $K(x)=0,75 \left(1-x^2 \right)\mathbf{1}_{[-1,\,1]}(x)$. Regarding the choice of stepsize $(\gamma_n)$ and bandwidth $(h_n)$, we selected   parameters that satisfy the assumptions used in this paper : $(h_n)= \left(n^{-0.2}\right)$ and $(\gamma_n) =  \left(n^{-1}\right)$. 

To obtain the results, particularly those concerning the calculation of $R^2\left(\widehat{\beta}_j\right)$, and to evaluate the computational times, we began with an initial size $n = 100$. Variations were then made by incrementally adding $p$   additional data points. For the recursive estimator, each of these $p$  data points was recursively included one by one, as depicted in (\ref{eqnfgga1}). in contrast, for the non-recursive estimator the  size of sample will be simply $n+p$. In the following section, we present the results obtained from the simulated data.

\subsection{Results with simulated data}

The simulations will be conducted using the following models:
\medskip

\noindent \textbf{Model 1:} $Y = X_1+X_2+X_3+X_4+\varepsilon $;

\medskip

\noindent\textbf{Model 2:} $Y = \left(X_1+X_2+X_3+X_4\right)^3+\varepsilon$.
\medskip

\begin{figure}[!h]
\centering
\leavevmode
\includegraphics[scale=0.4,bb=380 550 160 40]{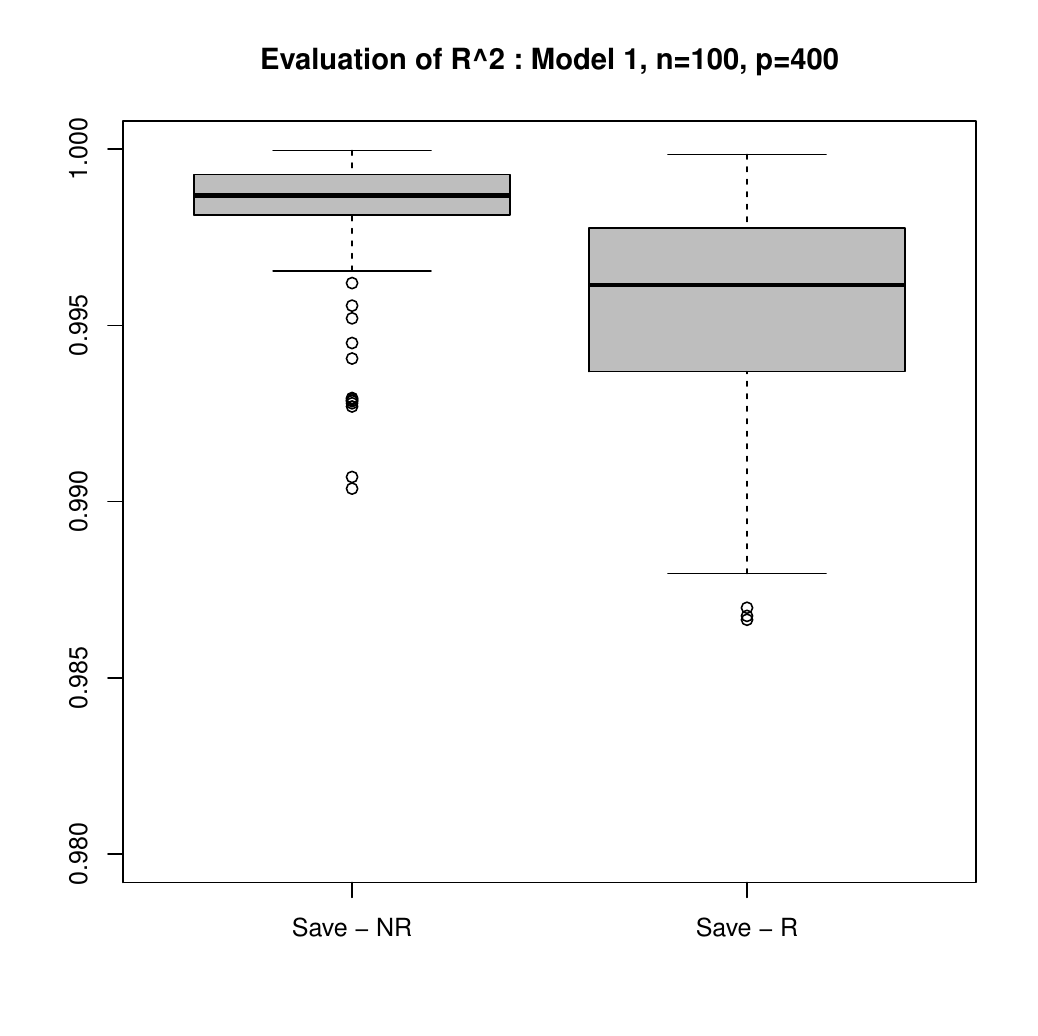}
\hspace{1cm}
\includegraphics[scale=0.4,bb=-250 550 160 40]{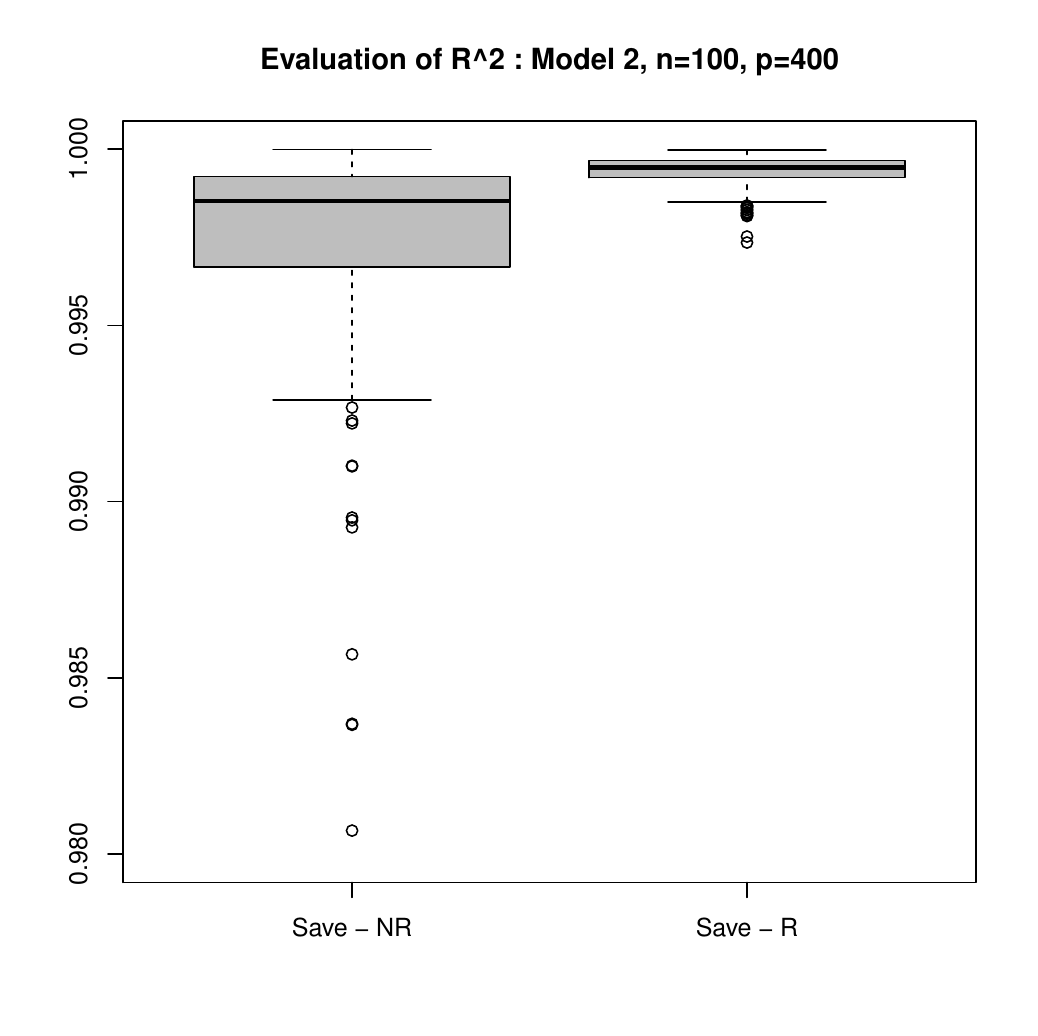}
\caption{Boxplots showing $R^2(\widehat{\beta}_1)$ for Model 1, $n = 100$, $p= 400$.}\label{box1}
\end{figure}

In both models $N=1$ and $\beta_1=(1,1,1,1,0)^T$. Each dataset was obtained as follows: $X=\left(X_1,X_2,X_3,X_4,X_5\right)^T$ is drawn from a multivariate normal distribution $\mathcal{N}(0,\textrm{\textbf{I}}_5)$, where $\textrm{\textbf{I}}_5$ is the $5\times 5$ identity matrix. $\varepsilon$ is generated from a standard normal distribution and $Y$ is  computed according to the  models described above. Table \ref{table1} below presents the means and standard deviations  of $R^2\left(\widehat{\beta}_1\right)$ obtained from the two models with  at least $200$ replications for $n=100$ and different values of $p$; $100$, $200$, $300$ and $400$. Figure 1 presents the  boxplots \ref{box1} of the estimators obtained for $n=100$ and $p=400$.

\begin{table}[h] 
\setlength{\tabcolsep}{0.2cm}  
\setlength{\doublerulesep}{0pt}
\renewcommand{\arraystretch}{0.6}  
\centering \caption{\label{table1} Means and standard deviations  of $R^2\left(\widehat{\beta}_1\right)$, over 200 replications for Models 1 and 2 with $n=100$, $p=100,\, 200,\, 300,\, 400$.}
\begin{center}
\scalebox{0.95}{
\begin{tabular}{p{1.5cm}p{1.5cm}lcccccccc}
\\\hline\hline\hline\\
        &                                      && &$p=100$  & &$p=200$ & &$p=300$& &$p= 400$        \\
\hline
\\                                                                                                         
\multirow{5}{*}{Model 1} & \multirow{5}{*}{$R^2\left(\widehat{\beta}_1\right)$}&$Save-NR$& & 0.96683  & &0.97291  & & 0.99676 & &0.99834 \\
        &                                                     &      & &(0.11677)& &(0.11602)&  & (0.00375)&   & (0.00157)\\
\\
				&                                                     &$Save-R$ & & 0.99178  & &0.99269  & & 0.99354 & &0.99554 \\
				&                                                     &      & &( 0.00460)& &(0.00371)&  & (0.00319)&   & (0.00291)\\
\\				
\hline
\\
\multirow{5}{*}{Model 2} & \multirow{5}{*}{$R^2\left(\widehat{\beta}_1\right)$}&$Save-NR$& &0.80008& &0.88389  & & 0.95250 & &0.96880 \\
        &                                                     &      & &(0.32103)& &( 0.26289)&  & ( 0.16695)&   & (0.12058) \\
\\
				&                                                     &$Save-R$ & & 0.99824  & &0.99895  & & 0.99905 & &0.99936   \\
				&                                                     &      & &( 0.00128)& &(0.00076)& & (0.00066)& & ( 0.00046)\\
\\
\hline				        
\hline\hline    
\end{tabular}
}
\end{center}
\end{table}

\vspace{6cm}

.  

\vspace{6cm}

We can already observe the superior performance of the recursive estimator ($Save-R$) over the non-recursive estimator  ($Save-NR$). Next, we evaluate the computational time taken by each estimator to perform the task.

\subsection{Comparison of computational times between non-recursive and recursive methods}

Let us clarify the execution of the estimators: the non-recursive estimator recalculates all computations each time one or more data points are added, whereas the recursive estimator adapts by leveraging previous computations. The computational times were evaluated based on this distinction. It is important to note that the results obtained are relative, as they depend on factors such as the computational software used and the specific computer hardware.

 Table \ref{table_times} presents the mean and standard deviation (in parentheses) of computational times for $n=100$ and for $p=25,\,p=50,\, p=75,\,p=100,\, 200,\, 300,\, 400$ based on $100$ replications. The results are valid for models 1 and 2. Note that the estimators that achieve relatively shorter computational times are those which insert the added data one by one. Therefore, the estimators of (\ref{eqnfgga1}) were used to obtain our results.
\begin{table}[h] 
\setlength{\tabcolsep}{0.2cm}  
\setlength{\doublerulesep}{0pt}
\renewcommand{\arraystretch}{1.1}  
\centering \caption{\label{table_times} Mean and standard deviation (in parentheses) computational times (in seconds) $\widehat{\beta}_{NR}$, $\widehat{\beta}_{R}$ for $n=100$(initialization) and added data for different values of $p$ obtained for 100 replications.}
\begin{center}
\scalebox{0.9}{
\begin{tabular}{p{0.5cm}lccccccccccccccc}
\\\hline\hline
$p=$                  & &$25$  & &$50$   & &$75$  & &$100$ & &$150$& &$200$& &$300$&&$400$ \\
\hline                                                                                                         
$\widehat{\beta}_{NR}$& & 31.35  & &43.88& &58.64 & &75.59 & &$117.62$& &$170.09$& &$334.52$&&500.04\\
                      & & (0.56)& & (0.69)& &(1.01)&&(0.59)&&(2.48) & & (3.45)& &  (28.45)& & (7.42)
\\
\hline 
$\widehat{\beta}_{R}$ & & 1.26  & &3.48  & &7.42   & &12.61 & &$28.40$& &$43.02$& &$100.46$&&175.93 \\
                      & & (0.09)& &(0.23)& & (0.45)& &(0.49)& &(0.36)& & (1.84)& & (1.39)&& (1.75)
\\
\hline\hline    
\end{tabular}
}
\end{center}
\end{table}

\subsection{Application to real data}
Here is an application using real data where we compare the two estimators. This dataset originates from the National Institute of Diabetes and Digestive and Kidney Diseases. The objective is to predict whether a patient has diabetes based on diagnostic measurements. The dataset consists of data from $2000$ female patients with diabetes. Below are the variables we consider: $Y= Glucose$ (Plasma glucose concentration over 2 hours in an oral glucose tolerance test)  as response variable and the predictor $X=\left(X_1, \cdots, X_6\right)$ such that $X_1=BMI$ (Body mass index), $X_2=Pregnancies$ (Number of times pregnant), $X_3=DPF$(Diabetes pedigree function : a function which scores likelihood of diabetes based on family history), $X_4=Age$ (years), $X_5=Insulin$, $X_6=BP$( Diastolic blood pressure). Therefore, we considered the situation where $Y=F\left(\beta^T X\right)$, and we estimated the direction $\beta = \left(\beta_1,\cdots,\beta_6\right)$ . Since we do not know the true direction to be estimated as in a simulation, the criterion used to assess the quality of the obtained estimators is as follows : we consider the vectors $BX_i = \left(\beta_1 X_{i1},\cdots,\beta_6X_{i6}\right)$ and $\widehat{B}X_i = \left(\widehat{\beta}_1 X_{i1},\cdots,\widehat{\beta}_6X_{i6}\right)$,. Thus we have
$$
R^2_i\left(\widehat{\beta}\right) = \frac{\left(\left(BX_i\right)^T \cdot\widehat{B}X_i\right)^2 }{\left(\left(\widehat{B}X_i\right)^T \cdot \widehat{B}X_i\right) \left(\left( BX_i\right)^T \cdot BX_i\right)}.
$$ 
It is  the squared multiple correlation coefficient between the variable $\left(BX_i\right)^T$ and the ideally reduced variables $\left(\widehat{B}X_i\right)^T$ (see \cite{li1991}).

 We proceeded as follows to obtain $\widehat{\beta}$: for the recursive estimator $(\widehat{\beta}_{R})$ we set $n=100$ and varied $p$ with different values. For the non-recursive estimator $(\widehat{\beta}_{NR})$ the  size of sample was $n+p$. We evaluated the quality of  $\widehat{\beta}$ by calculating $R^2_i\left(\widehat{\beta}_{R}\right)$ and $R^2_i\left(\widehat{\beta}_{NR}\right)$ on the remaining data, on  $2000 - (n+p)$ and $i\in \{n+p+1,\cdots, 2000 - (n+p)\}$. In  table (\ref{table2}), we present  the estimators $ \widehat{\beta}_{NR}$ and $ \widehat{\beta}_{R}$  obtained along with the means and standard errors (in parentheses) calculated on the remaining data.

\begin{table}[h] 
\setlength{\tabcolsep}{0.2cm}  
\setlength{\doublerulesep}{0pt}
\renewcommand{\arraystretch}{0.7}  
\centering \caption{\label{table2} $\widehat{\beta}_{NR}$, $\widehat{\beta}_{R}$, means and standard deviations  of $R_i^2\left(\widehat{\beta}\right)$, for $n=100$, $p=100,\, 200,\, 300,\, 400$.}
\begin{center}
\scalebox{0.9}{
\begin{tabular}{p{1cm}p{0.5cm}lccccccccccc||cc}
\\\hline\hline
                         &                     & &$\widehat{\beta}_1$  & &$\widehat{\beta}_2$ & &$\widehat{\beta}_3$& &$\widehat{\beta}_4$ & &$\widehat{\beta}_5$& &$\widehat{\beta}_6$& &Mean$(R^2_i)$ \\
\hline                                                                                                         
\multirow{2}{*}{p=100}  &$\widehat{\beta}_{NR}$& & 0.0849  & &-0.0067  & & 0.9956 & &-0.0136 & &$-0.0353$& &$0.0014$ & &$0.9843 \,(0.0590)$\\
				                &$\widehat{\beta}_{R}$ & & 0.1145  & &-0.0067  & & 0.9925 & &-0.0034 & &$-0.0386$& &$-0.0102$& &$0.9697\,(0.1046)$ \\
\hline                                                                                                         
\multirow{2}{*}{p=200}  &$\widehat{\beta}_{NR}$& & -0.0387  & &$-0.0229$  & & $-0.9989$ & &$0.0068$ & &$0.0134$& &$-0.0002$ & &$0.9726 \,(0.0993)$\\
				                &$\widehat{\beta}_{R}$ & &$-0.0433$& &$-0.0202$& &$-0.9988$& &$0.0036$& &$0.0125$& &$0.0010$& &$0.97262\,(0.0990)$ \\
\hline                                                                                                         
\multirow{2}{*}{p=300}  &$\widehat{\beta}_{NR}$& &$-0.0277$& &$-0.0083$& & $-0.9995$ & &$0.0070$& &$0.0117$& &$0.0001$ & &$0.9702 \,(0.1074)$\\
				                &$\widehat{\beta}_{R}$ & &$-0.0349$& &$-0.0047$& &$-0.9993$& &$0.0055$& &$0.0122$& &$0.0023$& &$0.9702\,(0.1074)$ \\
\hline                                                                                                         
\multirow{2}{*}{p=400}  &$\widehat{\beta}_{NR}$& &$0.0398$& &$-0.0021$& &$0.9990$& &$-0.0082$& &$-0.0141$& &$-0.0003$ & &$0.9704 \,(0.0997)$\\
				                &$\widehat{\beta}_{R}$ & &$0.0449$& &$-0.0036$& &$0.9989$& &$-0.0069$& &$-0.0139$& &$-0.0024$& &$0.9704\,(0.0994)$ \\	\hline\hline    
\end{tabular}
}
\end{center}
\end{table}

\bigskip

\section{Proofs}\label{proofs1}\label{section5}

In this section, we present the proof of Theorem \ref{theo1}. This proof requires preliminary results, including lemmas and propositions,  we will introduce at the very beginning of the section. Some of their proofs are deferred to the subsequent section.

\subsection{Preliminary results}
We begin with the three propositions below. Propositions \ref{lemma2a} and \ref{lemma3a} correspond to Propositions 6.2 and 6.3 in \cite{nkou2022c}, where they were formally established. Proposition \ref{lemma4a} is proven using a similar approach.

\begin{Proposition}\label{lemma2a}
Under Assumptions \textbf{A2.},  \textbf{A4.} and  \textbf{A5.}, we have
$$
\sup_{y\in\R}\left|f(y) - \widehat{f}_{n}(y)\right| = \mathbf{O}_{a.s.}\left(h_n^2 + \sqrt{\frac{\log n}{n h_n^2}}\right).
$$
\end{Proposition}

\begin{Proposition}\label{lemma3a}
Under Assumptions \textbf{A1.}, \textbf{A2.},  \textbf{A4.} and  \textbf{A5.}, we have for all $j\in \{1,\cdots,d\}$,
$$
\sup_{y\in\R}\left|\widehat{g}_{j,n}(y) - g_j(y) \right| = \mathbf{O}_{a.s.}\left(h_n^2 + M_n\sqrt{\frac{\log n}{nh_n^2}}\right).
$$
\end{Proposition}

\begin{Proposition}\label{lemma4a} 
Under Assumptions \textbf{A1.}, \textbf{A2.},  \textbf{A4.} and  \textbf{A5.}, we have for all $1\leq k,\ell \leq d$,
$$
\sup_{y\in\R}\left|\widehat{G}_{k\ell,n}(y) - G_{k\ell}(y) \right| = \mathbf{O}_{a.s.}\left(h_n^2 + M_n^2\sqrt{\frac{\log n}{nh_n^2}}\right).
$$
\end{Proposition}
We continue with the lemmas.

\begin{Lemma}\label{lemma01a}
We consider the sequence $\pi_n$ defined in Equation (\ref{eqnpi1}).It holds for all $n\in \N^*$,
$$
\pi_n\sum_{k=1}^n  \pi_k^{-1} \gamma_k=1-\pi_n.
$$
\end{Lemma}
\textbf{Proof.} Using proof by induction we show that \, $\sum_{k=1}^n  \pi_k^{-1} \gamma_k  + 1 = \pi_n^{-1}$.\, $\square$\\

\begin{Lemma}\label{lemma1a}
Let $(v_n) \in \mathcal{GS}(v^*)$, $(\gamma_n) \in \mathcal{GS}(-c)$, and let $m > 0$ such that $m - v^*\varepsilon_0 > 0$ where $\varepsilon_0$ is defined in (\ref{eqga1}). We have
$$
\lim_{n\rightarrow + \infty} v_n \left(\pi_n\right)^m \sum_{i=1}^n \left(\pi_i\right)^{-m} \frac{\gamma_i}{v_i} = \frac{1}{m - v^*\varepsilon_0}.
$$
Moreover, for all positive sequence $(\alpha_n)$ such that $\lim_{n\rightarrow +\infty} \alpha_n = 0$, and all $a\in \R$,
$$
\lim_{n\rightarrow + \infty} v_n \left(\pi_n\right)^m\left(\sum_{i=1}^n \left(\pi_i\right)^{-m} \frac{\gamma_i}{v_i} \alpha_i + a\right)=0.
$$
\end{Lemma}
Lemma \ref{lemma1a}  was established in \cite{mokkadem2009a}. It is widely used in stochastic approximation algorithms \cite{mokkadem2007, slaoui2015} and  is crucial throughout the proofs. 

\medskip

The proofs of the following lemmas are presented in the subsequent section.

\begin{Lemma}\label{rec_lema1}  Under Assumptions \textbf{A2.}, \textbf{A4.} and \textbf{A5.} we have
$$
\frac{b_n^{-2}}{\sqrt{n}}\sum_{i=1}^n \bigg(\widehat{H}_{1,n}(Y_i)  - H_1(Y_i)\bigg)\bigg(\widehat{H}_{2,n}(Y_i)  - H_2(Y_i)\bigg) = o_p(1),
$$
where $H_1(\cdot)$ and $H_2(\cdot)$   represent functions :$f(\cdot)$ , $g_j(\cdot)$ and $G_{k\ell}(\cdot)$ for $j\in \{1,\cdots,d\}$ and for each pair $1\leq k,\ell \leq d$.
\end{Lemma}
\textbf{Proof.} We will prove this lemma in the scenario where $H_1(\cdot)$ and $H_2(\cdot)$ are identical. Once this case is established, the other cases can be easily proven using Cauchy's inequality. Specifically, we will focus the case where $H_1(\cdot) = H_2(\cdot)= g_j(\cdot)$ for $j\in \{1,\cdots, d\}$. The proof for the other cases follows similarly. For the specified case , it is sufficient to demonstrate that :
$$
\frac{b_n^{-2}}{\sqrt{n}}\sum_{i=1}^n \bigg(\widehat{g}_{j,n}(Y_i)  - g_j(Y_i)\bigg)^2 = o_p(1).
$$

Proposition \ref{lemma3a}, under Assumption \textbf{A5.} guarantees that $
\frac{b_n^{-2}}{\sqrt{n}}\sum_{i=1}^n \bigg(\widehat{g}_{j,n}(Y_i)  - g_j(Y_i)\bigg)^2 = O_p\left(n^{-2c_2}\right)$,  which concludes the proof of the lemma.\, $\square$

\medskip

Below are three crucial lemmas. Their proofs are similar, so in  Section \ref{section6}, we will present only one them, specifically Lemma \ref{rec_lema4}.

\begin{Lemma}\label{rec_lema2} Under Assumptions \textbf{A.2}, \textbf{A.4} and \textbf{A.5}, we have :
$$
\frac{1}{\sqrt{n}}\sum_{i=1}^n \Big(H(Y_i)\widehat{f}_{n}(Y_i)\, -\, H(Y_i)f(Y_i)\Big) = \sqrt{n}\sum_{i=1}^n\Big(\pi_n\pi_i^{-1} \gamma_i H(Y_i) f(Y_i) \,-\, \frac{1}{n}\E(H(Y)f(Y))\Big)\,+\, o_p(1),
$$
where $H(\cdot)$ is one the following functions: $\frac{R_{k\ell}(\cdot)}{f_{b_n}(\cdot)}$, $\frac{r_k(\cdot)r_\ell(\cdot)}{f_{b_n}(\cdot)}$, $\frac{R_{kj}(\cdot)R_{j\ell}(\cdot) }{f_{b_n}(\cdot)}$,  $\frac{r_\ell(\cdot)r_k(\cdot)r_j(\cdot)}{f_{b_n}(\cdot)}$ and $\frac{R_{k\ell}(\cdot)r_k(\cdot)r_\ell(\cdot) }{f_{b_n}(\cdot)}$ with $(k,\, \ell, \,j) \in \{1,\cdots, d\}^3$.
\end{Lemma}

\medskip

\begin{Lemma}\label{rec_lema3} Under Assumptions \textbf{A.2},  \textbf{A.4} and \textbf{A.5}, we have:
$$
\frac{1}{\sqrt{n}}\sum_{i=1}^n \Big(H(Y_i)\widehat{g}_{k,n}(Y_i)\, -\, H(Y_i) g_k(Y_i)\Big) = \sqrt{n}\sum_{i=1}^n\Big(\pi_n\pi_i^{-1} \gamma_i X_{i,k} H(Y_i) f(Y_i) \,-\, \frac{1}{n}\E(H(Y)g_k(Y))\Big)\,+\, o_p(1),
$$
where $H(\cdot)$  is one the following functions: $\frac{r_\ell(\cdot)}{f_{b_n}(\cdot)}$, $\frac{r_\ell(\cdot)r_k(\cdot)}{f_{b_n}(\cdot)}$,  $\frac{r_\ell(\cdot)r_k(\cdot)r_j(\cdot)}{f_{b_n}(\cdot)}$ and $\frac{R_{kj}(\cdot)r_\ell(\cdot) }{f_{b_n}(\cdot)}$ with $(k,\, \ell, \,j) \in \{1,\cdots, d\}^3$.
\end{Lemma}

\medskip

\begin{Lemma}\label{rec_lema4} Under Assumptions \textbf{A.2},  \textbf{A.4} and \textbf{A.5}, we have :
\begin{multline*}
\frac{1}{\sqrt{n}}\sum_{i=1}^n \Big(H(Y_i)\widehat{G}_{k\ell,n}(Y_i)\, -\, H(Y_i) G_{k\ell}(Y_i)\Big)\\
 =\, \sqrt{n}\sum_{i=1}^n\Big(\pi_n\pi_i^{-1} \gamma_i X_{i,k}\, X_{i,\ell} H(Y_i) f(Y_i) \,-\, \frac{1}{n}\E(H(Y)G_{k\ell}(Y))\Big)\,+\, o_p(1),
\end{multline*}
where $H(\cdot)$ is one the following functions: $\frac{1}{f_{b_n}(\cdot)}$, $\frac{R_{k\ell}(\cdot)}{f_{b_n}(\cdot)}$,  and $\frac{r_k(\cdot)r_\ell(\cdot) }{f_{b_n}(\cdot)}$ with $(k,\, \ell) \in \{1,\cdots, d\}^2$.
\end{Lemma}

\subsection{Proof of Theorem \ref{theo1}}\label{section5_2}

In the following, $C$ stands a finite positive constant, whose specific value  is arbitrary and may vary between equations.

\medskip
For $n$ large enough, we have the following inequalities:
\begin{equation}\label{eqn1mn}
\begin{split}
\sup_{y\in \R}\left|g_j(y)\right|\leqslant  C M_n   \quad &  \quad  \sup_{y\in \R}\left|G_{k\ell}(y)\right|\leqslant C M_n^2,\\
\sup_{y\in \R}\left|\widehat{g}_{j,n}(y)\right|\leqslant  C M_n \quad &  \quad  \sup_{y\in \R}\left|\widehat{r}_{j,b_n}(y)\right|\leqslant C  b_n^{-1} M_n,\\
\quad \sup_{y\in \R}\left|\widehat{G}_{k\ell,n}(y)\right|\leqslant C M_n^2 \quad &  \quad \sup_{y\in \R}\left|\widehat{R}_{k\ell,b_n}(y)\right|\leqslant C b_n^{-1} M_n^2.
\end{split}
\end{equation}

These inequalities can be derived using the Assumptions \textbf{A1.}, \textbf{A2.}, \textbf{A4.}, \textbf{A5.}, following the method used in \cite{nkou2019} in the proof of Lemma 4.2.

We will develop the estimator $\widehat{\lambda}_{k\ell,n}$ defined in (\ref{eqnest1}).

For any $1\leq  k,\ell \leq  d$  and $j=1,\cdots,d$, we define the following functions:
\medskip

$
A_{k\ell,n}^{(1)}(y)=\frac{G_{k\ell}(y) \left(f_{b_n}(y)-\widehat{f}_{b_n}(y)\right)}{f_{b_n}(y)\widehat{f}_{b_n}(y)} + \frac{\widehat{G}_{k\ell,n} (y) - G_{k\ell}(y)}{f_{b_n}(y)} + \frac{G_{k\ell}(y)}{f_{b_n}(y)},\\
$

$
A_{k\ell,n}^{(2)}(y)=\frac{2g_\ell(y)g_k(y)\left(f_{b_n}(y)-\widehat{f}_{b_n}(y)\right)}{f_{b_n}(y)\widehat{f}_{b_n}(y)^2} + \frac{g_\ell(y)\bigg(\widehat{g}_{k,n}(y)-g_k(y)\bigg)}{f_{b_n}(y)^2}\\
$

$
		   \hspace{4cm} +\,  \frac{g_k(y)\bigg(\widehat{g}_{\ell,n}(y)-g_\ell(y)\bigg)}{f_{b_n}(y)^2} +\frac{g_\ell(y)g_k(y)\bigg)}{f_{b_n}(y)^2},\\
$

$			
B_{jk\ell,n}^{(1)}(y)= R_{j\ell,b_n}(y)\left(\widehat{R}_{k j,b_n}(y) -R_{kj,b_n}(y) \right) + R_{kj,b_n}(y) \left(\widehat{R}_{j \ell,b_n}(y)-R_{j\ell,b_n}(y)\right) +  R_{kj,b_n}(y) R_{j\ell,b_n}(y),\\
$

$
B_{jk\ell,n}^{(2)}(y) = \left(\widehat{R}_{k j,b_n}(y) - R_{kj,b_n}(y)\right)r_{j,b_n}(y)r_{\ell,b_n}(y)+ \big(\widehat{r}_{j,b_n}(y)- r_{j,b_n}(y)\big)R_{kj,b_n}(y)r_{\ell,b_n}(y)\\
$

$
																																			 \hspace{3cm} + \quad\bigg(\widehat{r}_{\ell,b_n}(y)- r_{\ell,b_n}(y)\bigg)r_{j,b_n}(y)R_{kj,b_n}(y) + r_{j,b_n}(y)r_{\ell,b_n}(y)R_{kj,b_n}(y)\\
$

$
B_{jk\ell,n}^{(3)}(y) = \left(\widehat{R}_{j \ell,b_n}(y) - R_{j\ell,b_n}(y)\right)r_{j,b_n}(y)r_{\ell,b_n}(y)+ \big(\widehat{r}_{j,b_n}(y)- r_{j,b_n}(y)\big)R_{kj,b_n}(y)r_{\ell,b_n}(y)\\
$

$
																																			 \hspace{3cm} + \quad\Bigg(\widehat{r}_{\ell,b_n}(y)- r_{\ell,b_n}(y)\Bigg)r_{j,b_n}(y)R_{j\ell,b_n}(y) + r_{j,b_n}(y)r_{\ell,b_n}(y)R_{j\ell,b_n}(y)
$

$
B_{jk\ell,n}^{(4)}(y)= r_{j,b_n}^2(y)r_{\ell,b_n}(y)\bigg(\widehat{r}_{k,b_n}(y)-r_{k,b_n}(y)\bigg) + r_{j,b_n}^2(y)r_{k,b_n}(y)\bigg(\widehat{r}_{\ell,b_n}(y)-r_{\ell,b_n}(y)\bigg)
$

$
\hspace{3cm} + \quad r_{j,b_n}^2(y)r_{k,b_n}(y)r_{\ell,b_n}(y) +2r_{j,b_n}(y)r_{k,b_n}(y)r_{\ell,b_n}(y)\bigg(\widehat{r}_{j,b_n}(y)-r_{j,b_n}(y)\bigg),
$

\medskip

The proof of the theorem will proceed in three steps.

\medskip

\textbf{Step 1.} The first step involves demonstrating that
\begin{multline}\label{eqstep1a}
\sqrt{n}\widehat{\lambda}_{k,\ell}^{(n)} = \frac{1}{\sqrt{n}}\sum_{i=1}^n\Bigg(\delta_{k,\ell}-2\left(A_{k\ell,n}^{(1)}(Y_i)-A_{k\ell,n}^{(2)}(Y_i) \right)\\
+\sum_{j=1}^d \left(B_{jk\ell,n}^{(1)}(Y_i) - B_{kj\ell,n}^{(2)}(Y_i)- B_{j\ell k}^{(2)}(Y_i) - B_{k j\ell}^{(3)}(Y_i)- B_{j\ell k}^{(3)}(Y_i) + B_{jk\ell,n}^{(4)}(Y_i)\right)\Bigg) + o_p(1).
\end{multline} 
Which is obtained using  Lemma \ref{rec_lema1} and inequalities (\ref{eqn1mn}). 

\medskip

\textbf{Step 2.} Show that
\begin{multline}\label{eqn1step2}
\frac{1}{\sqrt{n}}\sum_{i=1}^n\left(A_{k\ell,n}^{(1)}(Y_i) - \E\left(R_{k\ell,b_n}(Y)\right)\right)\\
 = \sqrt{n}\left(\pi_n \sum_{i=1}^n  \pi_i^{-1}\gamma_i \big(X_{i,k}X_{i,\ell}- R_{k\ell,b_n}(Y_i)\big)
 + \frac{1}{n}\sum_{i=1}^n\big(R_{k\ell,b_n}(Y_i) - \E\left(R_{k\ell,b_n}(Y)\right)\big)\right)+ o_p(1).
\end{multline}

Firstly, we observe that
\begin{multline*}
\sqrt{n}\left(\frac{G_{k\ell}(Y)\left(f_{b_n}(Y) - \widehat{f}_{b_n}(Y)\right)}{f_{b_n}(Y)\widehat{f}_{b_n}(Y)}\,-\,\frac{G_{k\ell}(Y)\left(f_{b_n}(Y) - \widehat{f}_{b_n}(Y)\right)}{f_{b_n}^2(Y)}\right)\\
 = \frac{\sqrt{n}G_{k\ell}(Y)\left(f_{b_n}(Y) - \widehat{f}_{b_n}(Y)\right)^2}{f_{b_n}^2(Y)\widehat{f}_{b_n}(Y)}.
\end{multline*}
However, according to Lemma \ref{lemma2a} , we obtain
$$
\frac{\sqrt{n}\left(f_{b_n}(Y) - \widehat{f}_{b_n}(Y)\right)^2}{f_{b_n}^2(Y)\widehat{f}_{b_n}(Y)} = O_p\left(n^{-1/2 + 2c_1 + 3c_2} \log n\right).
$$
According to Assumption \textbf{A5.} we have  $n^{-1/2 + 2c_1 + 3c_2} < n^{-c_2}$, thus
\begin{equation}\label{eqn2step2}
\frac{1}{\sqrt{n}}\sum_{i=1}^n\left(\frac{G_{k\ell}(Y_i)\left(f_{b_n}(Y_i)-\widehat{f}_{b_n}(Y_i)\right)}{f_{b_n}(Y_i)\widehat{f}_{b_n}(Y_i)}\right)\,=\,\frac{1}{\sqrt{n}}\sum_{i=1}^n\left(\frac{G_{k\ell}(Y_i)\left(f_{b_n}(Y_i) - \widehat{f}_{b_n}(Y_i)\right)}{f_{b_n}^2(Y_i)}\right)+o_p(1).
\end{equation}
Consequently
$$
\frac{1}{\sqrt{n}}\sum_{i=1}^n A_{k\ell,n}^{(1)}(Y_i) = \frac{1}{\sqrt{n}}\sum_{i=1}^n\left(\frac{\widehat{G}_{k\ell,n}(Y_i) - G_{k\ell}(Y_i)}{f_{b_n}(Y_i)}  + \frac{G_{k\ell}(Y_i)}{f_{b_n}(Y_i)} -  \frac{G_{k\ell}(Y_i)\left(\widehat{f}_{b_n}(Y_i) - f_{b_n}(Y_i)\right)}{f_{b_n}(Y_i)^2}\right) +o_p(1). 
$$
From Lemma \ref{rec_lema4}, using $H(Y)= \frac{1}{f_{b_n}(Y)}$, we derive
$$
\frac{1}{\sqrt{n}}\sum_{i=1}^n\left(\frac{\widehat{G}_{k\ell,n}(Y_i) - G_{k\ell}(Y_i)}{f_{b_n}(Y_i)}\right) = \sqrt{n}\sum_{i=1}^n  \left( \pi_n \pi_i^{-1}\gamma_i X_{i,k}X_{i,\ell} - \frac{1}{n}\E\left(R_{k\ell,b_n}(Y)\right)\right)+ o_p(1),
$$
and from Lemma \ref{rec_lema2} with $H(Y)= \frac{G_{k\ell}(Y)}{f_{b_n}^2(Y)}$, we have
$$
\frac{1}{\sqrt{n}}\sum_{i=1}^n\left(\frac{G_{k\ell}(Y_i)\left(\widehat{f}_{b_n}(Y_i) - f_{b_n}(Y_i)\right)}{f_{b_n}(Y_i)^2}\right) = \sqrt{n}\sum_{i=1}^n\left(\pi_n \pi_i^{-1}\gamma_i R_{k\ell,b_n}(Y_i) - \frac{1}{n}\E\left(R_{k\ell,b_n}(Y)\right)\right) + o_p(1),
$$
it follows that
$$
\frac{1}{\sqrt{n}}\sum_{i=1}^n A_{k\ell,n}^{(1)}(Y_i) = \sqrt{n}\sum_{i=1}^n \pi_n \pi_i^{-1}\gamma_i \left(X_{i,k}X_{i,\ell}- R_{k\ell,b_n}(Y_i)\right) + \frac{1}{\sqrt{n}}\sum_{i=1}^n R_{k\ell,b_n}(Y_i)+ o_p(1),
$$
and (\ref{eqn1step2}) is obtained.

Results (\ref{eqn1step3}) to (\ref{eqn1step6}) presented below are derived using the same method as for (\ref{eqn1step2}). We will not delve into the details here.

\medskip

\begin{multline}\label{eqn1step3}
\frac{1}{\sqrt{n}}\sum_{i=1}^n\left(A_{k\ell,n}^{(2)}(Y_i) - \E\left(r_{k,b_n}(Y)r_{\ell,b_n}(Y)\right)\right) \\
= \sqrt{n}\Bigg(\pi_n \sum_{i=1}^n \pi_i^{-1}\gamma_i\big(X_{i,k}r_{\ell,b_n}(Y_i) + X_{i,\ell}r_{k,b_n}(Y_i) - 2r_{k,b_n}(Y_i)r_{\ell,b_n}(Y_i)\big)\\
  +\, \frac{1}{n}\sum_{i=1}^n\big(r_{k,b_n}(Y_i)r_{\ell,b_n}(Y_i)-\E\left(r_{k,b_n}(Y)r_{\ell,b_n}(Y)\right)\big)\Bigg) + o_p(1);
\end{multline}

\begin{multline}\label{eqn1step4}
\frac{1}{\sqrt{n}}\sum_{i=1}^n\sum_{j=1}^d \left(B_{jk\ell,n}^{(1)}(Y_i) - \E\left(R_{kj,b_n}(Y_i) R_{j\ell,b_n}(Y_i)\right)\right)\\
 =\sqrt{n}\Bigg(\pi_n\sum_{i=1}^n \sum_{j=1}^d \pi_i^{-1}\gamma_i \Big(X_{i,\ell}X_{i,j}R_{j k,b_n}(Y_i) + X_{i,k}X_{i,j}R_{j\ell,b_n}(Y_i)- 2R_{j\ell,b_n}(Y_i)R_{kj,b_n}(Y_i)\Big)\\
+\, \frac{1}{n}\sum_{i=1}^n \sum_{j=1}^d\left(R_{j\ell,b_n}(Y_i)R_{kj,b_n}(Y_i)- \E\Big(R_{kj,b_n}(Y_i) R_{j\ell,b_n}(Y_i)\Big)\right)\Bigg)  + o_p(1);
\end{multline}

\begin{multline}\label{eqn1step5}
\frac{1}{\sqrt{n}}\sum_{i=1}^n\sum_{j=1}^d \left(B_{jk\ell,n}^{(2)}(Y_i) - \E\left(r_{\ell,b_n}(Y) R_{kj,b_n}(Y)r_{j,b_n}(Y)\right)\right)\\
= \sqrt{n}\Bigg(\pi_n \sum_{i=1}^n\sum_{j=1}^d \pi_i^{-1}\gamma_i\bigg(X_{i,k}X_{i,j}r_{j,b_n}(Y_i)r_{\ell,b_n}(Y_i) + X_{i,j}R_{kj,b_n}(Y_i)r_{\ell,b_n}(Y_i)\\
 + X_{\ell,i}R_{kj,b_n}(Y_i)r_{j,b_n}(Y_i) -3r_{\ell,b_n}(Y_i)R_{kj,b_n}(Y_i)r_{j,b_n}(Y_i)\bigg)\\
 + \frac{1}{n}\sum_{i=1}^n\sum_{j=1}^d \Big(r_{j,b_n}(Y_i)r_{\ell,b_n}(Y_i)R_{kj,b_n}(Y_i)\\
  - \E\left(r_{\ell,b_n}(Y) R_{kj,b_n}(Y)r_{j,b_n}(Y)\right)\Big)\Bigg) + o_p(1);
\end{multline}

\medskip

$B_{jk\ell,n}^{(3)}(Y)$ is obtained from $B_{jk\ell,n}^{(2)}(y)$ by permuting $k$ and $\ell$, thus
\begin{multline}\label{eqn1step5b}
\frac{1}{\sqrt{n}}\sum_{i=1}^n\sum_{j=1}^d \left(B_{jk\ell,n}^{(3)}(Y_i) - \E\left(r_{k,b_n}(Y) R_{\ell j,b_n}(Y)r_{j,b_n}(Y)\right)\right)\\
= \sqrt{n}\Bigg(\pi_n \sum_{i=1}^n\sum_{j=1}^d \pi_i^{-1}\gamma_i\bigg(X_{i,\ell}X_{i,j}r_{j,b_n}(Y_i)r_{k,b_n}(Y_i) + X_{i,j}R_{\ell j,b_n}(Y_i)r_{k,b_n}(Y_i)\\
 + X_{k,i}R_{\ell j,b_n}(Y_i)r_{j,b_n}(Y_i) -3r_{k,b_n}(Y_i)R_{\ell j,b_n}(Y_i)r_{j,b_n}(Y_i)\bigg)\\
 + \frac{1}{n}\sum_{i=1}^n\sum_{j=1}^d \Big(r_{j,b_n}(Y_i)r_{k,b_n}(Y_i)R_{\ell j,b_n}(Y_i)\\
  - \E\left(r_{k,b_n}(Y) R_{\ell j,b_n}(Y)r_{j,b_n}(Y)\right)\Big)\Bigg) + o_p(1);
\end{multline}

\medskip

\begin{multline}\label{eqn1step6}
\frac{1}{\sqrt{n}}\sum_{i=1}^n \sum_{j=1}^d \left(B_{jk\ell,n}^{(4)}(Y_i) - \E\left(r_{j,b_n}^2(Y)r_{k,b_n}(Y)r_{\ell,b_n}(Y)\right)\right)\\
= \sqrt{n}\Bigg(\pi_n\sum_{i=1}^n \pi_i^{-1}\gamma_i\bigg(X_{i,k}r_{j,b_n}^2(Y_i)r_{\ell,b_n}(Y_i) + X_{i,\ell}r_{j,b_n}^2(Y_i)r_{k,b_n}(Y_i)+2X_{i,j}r_{j,b_n}(Y_i)r_{k,b_n}(Y_i)r_{\ell,b_n}(Y_i)\\
 - 4r_{j,b_n}^2(Y_i)r_{k,b_n}(Y_i)r_{\ell,b_n}(Y_i)\bigg)\\
\, + \,\frac{1}{n}\sum_{i=1}^n \sum_{j=1}^d \bigg(r_{j,b_n}^2(Y_i)r_{k,b_n}(Y_i)r_{\ell,b_n}(Y_i)-\E\left(r_{j,b_n}^2(Y)r_{k,b_n}(Y)r_{\ell,b_n}(Y)\right)\bigg)\Bigg)\,+\, o_p(1).
\end{multline}

\medskip

\textbf{Step 3.} Conclusion of the proof of this theorem.\\
We define the following random variables:
\begin{eqnarray*}
\mathcal{A}^{(1)}_{k\ell,b_n}&=&\pi_n \sum_{i=1}^n  \pi_i^{-1}\gamma_i \big(X_{i,k}X_{i,\ell}- R_{k\ell,b_n}(Y_i)\big)\\
\mathcal{A}^{(2)}_{k\ell,b_n}&=& \pi_n \sum_{i=1}^n \pi_i^{-1}\gamma_i\big(X_{i,k}r_{\ell,b_n}(Y_i) + X_{i,\ell}r_{k,b_n}(Y_i) - 2r_{k,b_n}(Y_i)r_{\ell,b_n}(Y_i)\big)\\
\mathcal{B}^{(1)}_{jk\ell,b_n}&=&\pi_n\sum_{i=1}^n \pi_i^{-1}\gamma_i \Big(X_{i,\ell}X_{i,j}R_{k j,b_n}(Y_i) + X_{i,k}X_{i,j}R_{j\ell,b_n}(Y_i)- 2R_{j\ell,b_n}(Y_i)R_{kj,b_n}(Y_i)\Big)\\
 \mathcal{B}^{(2)}_{(j\ell) k,b_n}&=&\pi_n \sum_{i=1}^n \pi_i^{-1}\gamma_i\bigg(X_{i,\ell}X_{i,j}r_{j,b_n}(Y_i)r_{k,b_n}(Y_i) + X_{i,j}R_{j \ell,b_n}(Y_i)r_{k,b_n}(Y_i)\\
 & & \hspace{5cm} +\, X_{k,i}R_{\ell j,b_n}(Y_i)r_{j,b_n}(Y_i) -3r_{k,b_n}(Y_i)R_{j\ell,b_n}(Y_i)r_{j,b_n}(Y_i)\bigg)\\
\mathcal{B}^{(3)}_{(j\ell) k,b_n}&=&\pi_n \sum_{i=1}^n\sum_{j=1}^d \pi_i^{-1}\gamma_i\bigg(X_{i,\ell}X_{i,j}r_{j,b_n}(Y_i)r_{k,b_n}(Y_i) + X_{i,j}R_{\ell j,b_n}(Y_i)r_{k,b_n}(Y_i)\\
 & & \hspace{5cm} +\, X_{k,i}R_{j\ell,b_n}(Y_i)r_{j,b_n}(Y_i) -3r_{k,b_n}(Y_i)R_{j\ell ,b_n}(Y_i)r_{j,b_n}(Y_i)\bigg)\\
\mathcal{B}^{(4)}_{j k\ell,b_n}&=&\pi_n\sum_{i=1}^n \pi_i^{-1}\gamma_i\bigg(X_{i,k}r_{j,b_n}^2(Y_i)r_{\ell,b_n}(Y_i) + X_{i,\ell}r_{j,b_n}^2(Y_i)r_{k,b_n}(Y_i)+2X_{i,j}r_{j,b_n}(Y_i)r_{k,b_n}(Y_i)r_{\ell,b_n}(Y_i)\\
  & & \hspace{5cm} - \, 4r_{j,b_n}^2(Y_i)r_{k,b_n}(Y_i)r_{\ell,b_n}(Y_i)\bigg).
\end{eqnarray*}
Considering (\ref{eqstep1a}) and the random variables below, we define the $k\ell$th element of matrix $\Gamma_{b_n}$ as
$$
\Gamma_{k\ell,b_n}=\delta_{k,\ell}-2\left(\mathcal{A}_{k\ell,b_n}^{(1)}-\mathcal{A}_{k\ell,b_n}^{(2)} \right)\\
+\sum_{j=1}^d\left(\mathcal{B}_{jk\ell,b_n}^{(1)} - \mathcal{B}_{(j\ell) k,b_n}^{(2)}- \mathcal{B}_{(kj) \ell,b_n}^{(2)} - \mathcal{B}_{(j\ell) k,b_n}^{(3)}- \mathcal{B}_{(k j) \ell,b_n}^{(3)} + \mathcal{B}_{jk\ell,b_n}^{(4)}\right),
$$
and  the $k\ell$th element of the matrix $T_{b_n}(X_i,Y_i)$ as
\begin{eqnarray*}
T_{k\ell,b_n}(X_i,Y_i)&=& \delta_{k\ell} -2 \big(X_{i,k}X_{i,\ell}+X_{i,k}r_{\ell,b_n}(Y_i)+X_{i,\ell}r_{k,b_n}(Y_i)- r_{k,b_n}(Y_i)r_{\ell,b_n}(Y_i)\big)\\
                  &&+ \sum_{j=1}^d\big(X_{i,k}X_{i,j}R_{j\ell,b_n}(Y_i) + X_{i,\ell}X_{i,j}R_{jk,b_n}(Y_i) - R_{kj,b_n}(Y_i)R_{j\ell,b_n}(Y_i)\big)\\
									&&+\sum_{j=1}^d\big(X_{i,k}X_{i,j}r_{j,b_n}(Y_i)r_{\ell,b_n}(Y_i) + X_{i,j}r_{\ell,b_n}(Y_i)R_{kj,b_n}(Y_i) + X_{i,\ell}r_{j,b_n}(Y_i)R_{kj,b_n}(Y_i)\\
									&& \hspace{8cm} -\, 2r_{\ell,b_n}(Y_i)r_{j,b_n}(Y_i)R_{kj,b_n}(Y_i) \big)\\
									&&+\sum_{j=1}^d\big(X_{i,j}X_{i,\ell} r_{\ell,b_n}(Y_i) + X_{i,j}r_{\ell,b_n}(Y_i)R_{kj,b_n}(Y_i) + X_{i,\ell}r_{j,b_n}(Y_i)R_{kj,b_n}(Y_i)\\
									&& \hspace{8cm}-\, 2r_{j,b_n}(Y_i)r_{\ell,b_n}(Y_i)R_{kj,b_n}(Y_i) \big)\\
									&&+\sum_{j=1}^d\big(X_{i,k}r_{j,b_n}^2(Y_i)r_{\ell,b_n}(Y_i) + X_{i,\ell}r_{j,b_n}^2(Y_i)r_{j,b_n}(Y_i)+2X_{i,j}r_{j,b_n}(Y_i)r_{k,b_n}(Y_i)r_{\ell,b_n}(Y_i)\\
									&& \hspace{8cm}-\, 3r_{j,b_n}^2(Y_i)r_{k,b_n}(Y_i)r_{\ell,b_n}(Y_i)\big).
\end{eqnarray*}
From (\ref{lambda_bn1}),  it is straightforward to observe that $\lambda_{k\ell,b_n} = \E\left(T_{k\ell,b_n}(X,Y)\right)$. \\
By defining $T_{k\ell,b_n}=\frac{1}{n}\sum_{i=1}^n T_{k\ell,b_n}(X_i,Y_i)$, from Steps 1 and 2, we obtain
$$
\sqrt{n}\left(\widehat{\lambda}_{k\ell,n} - \lambda_{k\ell,b_n}\right) = \sqrt{n}\Big(\Gamma_{k\ell,b_n} + T_{k\ell,b_n} - \E\left(\Gamma_{k\ell,b_n}+T_{k\ell,b_n}\right)\Big)\,+\, o_p(1).
$$ 
Clearly, 
\begin{equation}\label{eqncon1}
\sqrt{n}\lambda_{k\ell,b_n} = \sqrt{n}\lambda_{k\ell}\,+\, o_p(1)
\end{equation}
 and
\begin{multline}\label{eqncon2}
\sqrt{n}\Big(\Gamma_{k\ell,b_n} + T_{k\ell,b_n} - \E\big(\Gamma_{k\ell,b_n} + T_{k\ell,b_n}\big)\Big) = \sqrt{n}\Big(\Gamma_{k\ell,n} + T_{k\ell,n} - \E\big(\Gamma_{k\ell,n} + T_{k\ell,n}\big)\Big)\,+\, o_p(1).
\end{multline}
$\Gamma_{k\ell}$ and $T_{k\ell}$ denote the corresponding expressions of $\Gamma_{k\ell,b_n}$ and $T_{k\ell,b_n}$ respectively, obtained by substituting   $f_{b_n}$ with $f$. (\ref{eqncon1}) and  (\ref{eqncon2}) are derived straightforwardly using Assumptions \textbf{A3.} and  \textbf{A6.} and the arguments used by \cite{zhu1996}. For instance, refer to Equations (4.6) and (4.7) of the cited article.\\
Putting $\Phi_{n}=\sqrt{n}\left(\widehat{\Lambda}_n-\Lambda\right)$ and $\Phi_{k\ell,n}=\sqrt{n}\left(\widehat{\lambda}_{k\ell,n} - \lambda_{k\ell}\right)$, we have
$$
\textrm{tr} \Big(A^T \Phi_{n}\Big) = \sum_{k=1}^d\sum_{\ell =1}^d a_{k\ell}\Phi_{k\ell,n} =\sqrt{n}\sum_{k=1}^d\sum_{\ell =1}^d a_{k\ell}\Big(\Gamma_{k\ell,n} + T_{k\ell,n} - \E\left(\Gamma_{k\ell,n}+T_{k\ell,n}\right)\Big) \,+\, o_p(1).  
$$
Thus, by Central Limit Theorem and Slutsky's theorem, we deduce that  $\textrm{tr}\Big(A^T \Phi_{n}\Big)\stackrel{\mathscr{D}}{\longrightarrow} \mathcal{N}(0,\sigma^2_{A}) $ as $n\rightarrow +\infty$, where $\sigma^2_{A}$ is given in (\ref{siga}). Then, using Levy's theorem, we conclude that $\Phi_{n} \stackrel{\mathscr{D}}{\longrightarrow}  \Phi$ as $n\rightarrow +\infty$, $\Phi$ has a normal distribution in $\mathscr{M}_d(\R)$ with   $\textrm{tr}\Big(A^T \Phi\Big) \leadsto \mathcal{N}(0,\sigma^2_{A})$\, $\square$

\section{Proof of Lemma \ref{rec_lema4}}\label{section6}

The proofs of Lemmas \ref{rec_lema2} to \ref{rec_lema4} are essentially identucal,the proofs of the other lemmas can be easily deduced from this one.

\medskip

To establish their results, \cite{zhu2007} used a $U-$statistic sum. However, we do not have the same conditions required to employ this method, namely a sum of random variables with identical expectations. The proof will proceed in several steps.\\
We define the random variables
\begin{multline*}
w\left(X_{i,k},\, X_{i,\ell},Y_i,X_{j,k},\, X_{j,\ell},Y_j\right)\\
 =    H(Y_i) \pi_j^{-1} \gamma_j X_{j,k}\, X_{j,\ell} \frac{1}{h_j}K\left(\frac{Y_i-Y_j}{h_j}\right) + H(Y_j) \pi_i^{-1} \gamma_i X_{i,k} \, X_{i,\ell}\frac{1}{h_i}K\left(\frac{Y_j-Y_i}{h_i}\right) 
\end{multline*}
and \,
$
U_{n} = \frac{\pi_n}{n} \sum_{1\leq i < j \leq n} w\left(X_{i,k},\, X_{i,\ell},Y_i,X_{j,k},\, X_{j,\ell},Y_j\right).
$
\medskip 

\textbf{Step 1.} :Find a useful approximation of $\frac{1}{\sqrt{n}}\sum_{i=1}^n H(Y_i)\widehat{G}_{k\ell,n}(Y_i)$.\\
We have
$$
\frac{1}{\sqrt{n}}\sum_{i=1}^n H(Y_i)\widehat{G}_{k\ell,n}(Y_i) = \sqrt{n}U_{n} + \sqrt{n} \pi_n  U_{n,1},
$$
where\, $U_{n,1} = \frac{1}{n}\sum_{i=1}^n H(Y_i)X_{i,k}\, X_{i,\ell} \pi_i^{-1} \gamma_i \frac{1}{h_i}K\left(0\right)$. We have
$$
\E\left(\pi_n \sqrt{n} \, U_{n,1}\right) = K\left(0\right) \E\left(H(Y_1)X_{1,k}\, X_{1,\ell}\right) \left(\frac{\pi_n}{\sqrt{n}}\sum_{i=1}^n  \pi_i^{-1} \gamma_i \frac{1}{h_i}\right).
$$
From this equality, using Assumption \textbf{A1.} and Lemma \ref{lemma1a}, we obtain $ \left|\E\left(\pi_n \sqrt{n} \, U_{n,1}\right)\right| \sim M_n^2/\sqrt{n}h_n$. Therefore \, $\lim_{n\rightarrow + \infty}\E\left(\pi_n \sqrt{n} \, U_{n,1}\right) = 0$.\\
 Furthermore
$$
var\left(\pi_n \sqrt{n} \, U_{n,1}\right) = K^2\left(0\right)var\left(H(Y_i)X_{i,k}\, X_{i,\ell} \right)\left(\frac{\pi_n^2}{n}\sum_{i=1}^n \pi_i^{-2} \gamma_i^2 \frac{1}{h_i^2}\right).
$$ 
Still according to Lemma \ref{lemma1a}, we have $\left(\frac{\pi_n^2}{n}\sum_{i=1}^n \pi_i^{-2} \gamma_i^2 \frac{1}{h_i^2}\right) \sim \frac{\gamma_n}{nh_n^2}$, consequently \,$\lim_{n\rightarrow +\infty }var\left(\pi_n \sqrt{n} \, U_{n,1}\right) = 0$.
We can therefore deduce that \,$\pi_n \sqrt{n} \, U_{n,1} = o_p(1)$.\\
Clearly
\begin{equation}\label{eqnc1lem3}
\frac{1}{\sqrt{n}}\sum_{i=1}^n H(Y_i)\widehat{G}_{k\ell,n}(Y_i) = \sqrt{n}U_{n} + o_p(1).
\end{equation}
In the following, we will use index 1 in place of index $i$ and index 2 in place of index $j$ simply to indicate the scenario where $Y_i$ and $Y_j$ take different values. In this case, $h_1$, $\gamma_1$ and $\pi_1$ or alternatively $h_2$, $\gamma_2$ and $\pi_2$ will be taken as more indication than values resulting from the terms of the sequences $(h_n)$, $(\gamma_n)$ and $(\pi_n)$.

\medskip

Now we will develop $\E(U_{n})$. First, we have
\begin{eqnarray*}
& &\E\left[H(Y_1) \pi_2^{-1} \gamma_2 X_{2,k}\,X_{2,\ell} \frac{1}{h_2}K\left(\frac{Y_1-Y_2}{h_2}\right) \right]\\
& & \qquad = \pi_2^{-1} \gamma_2 \int \left(\int  K(u)\left(G_{k\ell}(y_1-h_2u) - G_{k\ell}(y_1)\right)du\right)H(y_1) f(y_1) dy_1\, + \, \pi_2^{-1} \gamma_2\E(H(Y)G_{k\ell}(Y)).
\end{eqnarray*}
Let $\sigma_3(h)=\int  K(u)\left(G_{k\ell}(y-h u) - G_{k\ell}(y)\right)du$, we obtain
$$
\E\left[H(Y_i) \pi_2^{-1} \gamma_2 X_{j,k}\,X_{j,\ell} \frac{1}{h_j}K\left(\frac{Y_i-Y_j}{h_j}\right) \right] =\pi_j^{-1} \gamma_j\E(H(Y)G_{k\ell}(Y)) + \E\left(H(Y)\right)\pi_j^{-1} \gamma_2\sigma_3(h_j).
$$
Using the same reasoning as for $\sigma_1$ mentioned above, we have $\sigma_3(h) \sim h^4$. Therefore,
$$
\E\left[w\left(X_{i,k},X_{i,\ell},Y_i,X_{j,k},X_{j,\ell},Y_j\right)\right] = (\pi_i^{-1} \gamma_i+\pi_j^{-1} \gamma_j)\E(H(Y)G_{k\ell}((Y)) + \E(H(Y))(\pi_i^{-1} \gamma_i \sigma_3(h_i)+\pi_j^{-1} \gamma_j \sigma_3(h_j)).
$$
And we have
\begin{eqnarray*}
\E(U_{n}) &=&\E(H(Y)G_{k\ell}(Y))\frac{\pi_n}{n} \sum_{1\leq i < j \leq n}\left(\pi_i^{-1} \gamma_i + \pi_j^{-1} \gamma_j\right) + \E(H(Y)) \frac{\pi_n}{n} \sum_{1\leq i < j \leq n}\left(\pi_i^{-1} \gamma_i \sigma_3(h_i) + \pi_j^{-1}\gamma_j \sigma_3(h_j)\right)\\
					&=&\E(H(Y)G_{k\ell}(Y)) \left(1-\frac{1}{n}\right)(1-\pi_n) + \E(H(Y)) \sigma_3(h_n) \left(1-\frac{1}{n}\right)\\
					&=&\E(H(Y)G_{k\ell}(Y)) + O\left(h_n^4 + \pi_n \right).
\end{eqnarray*}
As discussed earlier, we have $\pi_n \sim \gamma_n h_n^{-1}$ and from Assumption \textbf{A4.}, we derive $\gamma_n h_n^{-1} \sim n^{-c_3}$ where $3/4 < c_3 < 4/5$. Thus $h_n^{-4} \pi_n \sim  n^{-(c_3 - 4c_1)}$. Since $c_3 - 4c_1> 1/60$, it follows that
\begin{equation}\label{eqnlr5a}
\E(U_{n}) = \E(H(Y)G_{k\ell}(Y)) + O\left(h_n^4\right).
\end{equation}

\medskip 

\textbf{Step 2.}: We show that $U_{n}$ can be approximated by its projection.\\
Let
\begin{equation*}
\widehat{U}_n =\sum_{i=1}^n \bigg[\E\big(U_{n} |  \left(X_{i,k}, X_{i,\ell}, Y_i\right)\big)\bigg] \,- \,(n-1)\E\left(U_{n}\right),
\end{equation*}
and
$$
\widehat{U}_n(x_{k},x_{\ell},y) =\sum_{i=1}^n \bigg[\E\big(U_{n} |  \left(X_{i,k}, X_{i,\ell}, Y_i\right) = (x_{ik},x_{i\ell},y_i)\big)\bigg] \,- \,(n-1)\E\left(U_{n}\right).
$$
 We denote the following conditional expectations: 
$$U_1(x_{ik},x_{i\ell},y_i)= \E\big(U_{n} |  \left(X_{i,k}, X_{i,\ell} , Y_i\right) = (x_{ik},x_{i\ell},y_i)\big).$$
We have
\begin{eqnarray*}
& &U_1(x_{1k},x_{1\ell},y_1)\\
& & \qquad = \frac{\pi_n}{n}\sum_{i=1,2\leq j\leq n}\E\bigg[H(y_1) \pi_j^{-1} \gamma_j X_{j,k}\, X_{j,\ell} \frac{1}{h_j}K\left(\frac{y_1-Y_j}{h_j}\right)\\
						& &\hspace{1cm}   +\, H(Y_j) \pi_1^{-1} \gamma_1 x_{1k}x_{1\ell} \frac{1}{h_1}K\left(\frac{Y_j-y_1}{h_1}\right) \,|\,\left(X_{1,k},X_{1,\ell}, Y_1\right) = (x_{1k},x_{1\ell},y_1) \bigg]\\
						& & \hspace{1cm} +\,\frac{\pi_n}{n}\sum_{2\leq i < j \leq n}\E\bigg[  H(Y_i) \pi_j^{-1} \gamma_j X_{j,k}\, X_{j,\ell} \frac{1}{h_j}K\left(\frac{Y_i-Y_j}{h_j}\right) + H(Y_j) \pi_i^{-1} \gamma_i X_{i,k}\, X_{i,\ell} \frac{1}{h_i}K\left(\frac{Y_j-Y_i}{h_i}\right)\bigg]\\
						& & \qquad = U_{1a}(x_{1k},x_{1\ell},y_1)+U_{1b}x_{1k},x_{1\ell},y_1)+ U_{1c}(x_{1k},x_{1\ell},y_1).
\end{eqnarray*}
In the sequel, we will  detail $U_{1a}(x_{1k},x_{1\ell},y_1)$, $U_{1b}(x_{1k},x_{1\ell},y_1)$, and $U_{1c}(x_{1k},x_{1\ell},y_1)$. First, we have
$$
\E\bigg[H(y_1) \pi_j^{-1} \gamma_j X_{j,k} \, X_{j,\ell} \frac{1}{h_j}K\left(\frac{y_1-Y_j}{h_j}\right)\bigg] = \pi_j^{-1} \gamma_j H(y_1) G_{k\ell}(y_1) \,+ \, H(y_1) \pi_j^{-1} \gamma_j \sigma_3(h_j),
$$
and thus
\begin{eqnarray*}
U_{1a}(x_{1k},x_{1\ell},y_1) &=& \left(\frac{\pi_n}{n}\sum_{j=1}^n\pi_j^{-1}\gamma_j \right)H(y_1) G_{k\ell}(y_1)-H(y_1) G_{k\ell}(y_1)\frac{\pi_n}{n}\pi_1^{-1}\gamma_1 \\
			& & \qquad + \, H(y_1)\frac{\pi_n}{n}\sum_{j=1}^n \pi_j^{-1} \gamma_j \sigma_3(h_j) - H(y_1)\frac{\pi_n}{n}\pi_1^{-1} \gamma_1 \sigma_3(h_1).\\
&=&\left(\frac{\pi_n}{n}\sum_{j=1}^n\pi_j^{-1}\gamma_j \right)H(y_1) G_{k\ell}(y_1) + \frac{\pi_n}{n} + \frac{\sigma_3(h_n)}{n}\\
&=&\frac{1}{n}H(y_1) G_{k\ell}(y_1) + O\left(\frac{h_n^4}{n}\right).
\end{eqnarray*}
Similarly, we have
\begin{multline*}
\E\bigg[H(Y_j) \pi_1^{-1} \gamma_1 x_{1k}x_{1\ell}  \frac{1}{h_1} K\left(\frac{Y_j-y_1}{h_1}\right) \,|\,\left(X_{1,k}, X_{1,\ell}, Y_1\right) = (x_{1k},x_{1\ell} ,y_1)  \bigg]\\
 = \pi_1^{-1} \gamma_1 x_{1k} x_{1\ell}  H(y_1) f(y_1) + \pi_1^{-1} \gamma_1 x_{1k}x_{1\ell}  \sigma_3(h_1),
\end{multline*}
thus
\begin{eqnarray*}
U_{1b}(x_{1k},x_{1\ell},y_1) &=& \pi_n\left(1-\frac{1}{n}\right)\pi_1^{-1} \gamma_1 x_{1k}x_{1\ell} H(y_1) f(y_1)+ \pi_n\left(1-\frac{1}{n}\right)\pi_1^{-1} \gamma_1 x_{1k}x_{1\ell} \sigma_3(h_1)\\
			 &=& \pi_n\pi_1^{-1} \gamma_1 x_{1k}x_{1\ell} H(y_1) f(y_1) + \pi_n\pi_1^{-1} \gamma_1 x_{1k}x_{1\ell} \sigma_3(h_1)\\
			& & \qquad \qquad + \, \frac{\pi_n}{n}\left(\pi_1^{-1} \gamma_1 x_{1k}x_{1\ell} H(y_1) f(y_1) + \pi_1^{-1} \gamma_1 x_{1k}x_{1\ell} h_1^4\right)\\
			&=& \pi_n\pi_1^{-1} \gamma_1 x_{1k}x_{1\ell} H(y_1) f(y_1) + \pi_n\pi_1^{-1} \gamma_1 x_{1} \sigma_3(h_1) + \frac{\pi_n}{n}\pi_1^{-1} \gamma_1(1+ \sigma_3(h_1)).
\end{eqnarray*}
In the sequel
\begin{eqnarray*}
& & U_{1c}(x_{ik}, x_{i\ell},y_i)\\
& & \qquad =  \E\left(U_n\right) \,-\, \frac{\pi_n}{n}\sum_{j=2}^n\E\bigg[  H(Y_1) \pi_j^{-1} \gamma_j X_{j,k} X_{j,\ell}\frac{1}{h_j}K\left(\frac{Y_1-Y_j}{h_j}\right) + H(Y_j) \pi_1^{-1} \gamma_1 X_{1,k}X_{1,\ell} \frac{1}{h_1}K\left(\frac{Y_j-Y_1}{h_1}\right)\bigg]\\
& & \qquad = \E\left(U_n\right) \,-\, \frac{1-\pi_n}{n}\E(H(Y)G_{k\ell}(Y)) - \frac{\sigma_3(h_n)}{n} - \pi_n\pi_1^{-1} \gamma_1\E(H(Y)G_{k\ell}(Y))- \pi_n\E(H(Y))\pi_1^{-1} \gamma_1 \sigma_3(h_1)\\
& & \qquad = \E\left(U_n\right) \,-\,\frac{1}{n}\E(H(Y)G_{k\ell}(Y))- \pi_n\pi_1^{-1} \gamma_1\E(H(Y)G_{k\ell}(Y))\\
& & \hspace{3cm} + \frac{\pi_n}{n}\E(H(Y)G_{k\ell}(Y))- \pi_n\E(H(Y))\pi_1^{-1} \gamma_1 \sigma_3(h_1) - \frac{\sigma_3(h_n)}{n}.
\end{eqnarray*}
Therefore, we have
\begin{eqnarray*}
\sum_{i=1}^n U_{1c}(x_{ik},x_{i\ell},y_i)&=& n\E\left(U_n\right) - \E(H(Y)G_{k\ell}(Y)) - \left(\pi_n\sum_{i=1}^n \pi_i^{-1} \gamma_i\right)\E(H(Y)G_{k\ell}(Y))\\
& & \hspace{1cm} +\pi_n\E(H(Y)G_{k\ell}(Y)) - \left(\pi_n\sum_{i=1}^n\pi_i^{-1} \gamma_i \sigma_3(h_i)\right)\E(H(Y)) - \sigma_3(h_n)\\
                           &=& n\E\left(U_n\right) - 2\E(H(Y)G_{k\ell}(Y))+2\pi_n\E(H(Y)G_{k\ell}(Y)) - 2 \E(H(Y))\sigma_3(h_n)\\
													 &=& n\E\left(U_n\right) - 2\E(H(Y)G_{k\ell}(Y)) + \sigma_3(h_n).
\end{eqnarray*}

From the above results, we obtain
\begin{eqnarray*}
& & \widehat{U}_n(x_{k},x_{\ell},y) - \E\left(U_{n}\right)\\
& & \qquad = \sum_{i=1}^n U_1(x_{ik},x_{i\ell},y_i) - (n-1)\E\left(U_{n}\right) - \E\left(U_{n}\right) \\
& & \qquad = \sum_{i=1}^n \bigg(U_{1a}(x_{ik},x_{i\ell},y_i) + U_{1b}(x_{ik},x_{i\ell},y_i) + U_{1c}(x_{ik},x_{i\ell},y_i)\bigg)  - (n-1)\E\left(U_{n}\right) - \E\left(U_{n}\right)\\
& & \qquad = \sum_{i=1}^n \left[\frac{1}{n}H(y_i) G_{k\ell}(y_i) + \pi_n\pi_i^{-1} \gamma_i x_{ik}x_{i\ell} H(y_i) f(y_i)\right]- 2\E(H(Y)G_{k\ell}(Y)) + \sigma_3(h_n).
\end{eqnarray*}
A simplified form of $\widehat{U}_n$ is thus:
\begin{equation}\label{eqn1u}
\widehat{U}_n = \sum_{i=1}^n \left[\frac{1}{n}H(Y_i) G_{k\ell}(Y_i) + \pi_n\pi_i^{-1} \gamma_i X_{i,k} X_{i,\ell} H(Y_i) f(Y_i)\right]- \E(H(Y)G_{k\ell}(Y)) + O(h_n^4).
\end{equation}

\medskip

\textbf{Step 3.} : Express $U_n - \widehat{U}_n$ in a form similar to that of a   $U-$statistics writing.\\
We will consider the following conditional Expectation
$$
w_1(X_{t,k},X_{t,\ell},Y_t) = \E\left( w\left(X_{i,k}, X_{i,\ell},Y_i,X_{j,k},X_{j,\ell} ,Y_j\right) \,| \, (X_{t,k},X_{t,\ell},Y_t)\right),\qquad t=i,j.
$$
We have
\begin{eqnarray*}
U_n - \widehat{U}_n &=& U_n -\sum_{i_1=1}^n \bigg[\E\big(U_{n} |  \left(X_{i_1,k},X_{i_1,\ell}, Y_{i_1}\right)\big)\bigg] \,+ \,(n-1)\E\left(U_{n}\right)\\
										&=& \frac{\pi_n}{n} \sum_{1\leq i < j \leq n} w\left(X_{i,k},X_{i,\ell},Y_i,X_{j,k},X_{j,\ell},Y_j\right)-\frac{\pi_n}{2n}\sum_{1\leq i \neq j \leq n}w_1(X_{i,k},X_{i,\ell},Y_i)\\
										& &-\frac{\pi_n}{2n}\sum_{1\leq i \neq j \leq n}w_1(X_{j,k},X_{j,\ell},Y_j)\\
										& &\qquad \qquad -\frac{\pi_n}{2n}\sum_{(1\leq i \neq j \leq n)}\sum_{(i_1\neq i, i_1\neq j, 1\leq i_1\leq n)}\E\left( w\left(X_{i_1,k},X_{i_1,\ell},Y_{i_1},X_{j_1,k},X_{j_1,\ell},Y_{j_1}\right)\right) \,+\,(n-1)\E\left(U_{n}\right)\\
										&=& \frac{\pi_n}{n} \sum_{1\leq i < j \leq n} \left[w\left(X_{i,k},X_{i,\ell},Y_i,X_{j,k},X_{j,\ell},Y_j\right) - w_1(X_{i,k},X_{i,\ell},Y_i) - w_1(X_{j,k},X_{j,\ell},Y_j)\right]\\
										& &  -\sum_{(i_1\neq i, i_1\neq j, 1\leq i_1\leq n)} \frac{\pi_n}{n}\sum_{1\leq i  <  j \leq n}\E\left( w\left(X_{i,k},X_{i,\ell},Y_i,X_{j,k},Y_j\right)\right)\,+\,(n-1)\E\left(U_{n}\right)\\
										& & \qquad -\sum_{(i_1\neq i, i_1\neq j, 1\leq i_1\leq n)} \E\left(U_{n}\right)\,+\,(n-1)\E\left(U_{n}\right)\\
										&=&\frac{\pi_n}{n} \sum_{1\leq i < j \leq n} \left[w\left(X_{i,k},X_{i,\ell},Y_i,X_{j,k},X_{j,\ell},Y_j\right) - w_1(X_{i,k},X_{i,\ell},Y_i) - w_1(X_{j,k},X_{j,\ell},Y_j)\right]\\
										& & \qquad -(n-2) \E\left(U_{n}\right)\,+\,(n-1)\E\left(U_{n}\right)\\
										&=&\frac{\pi_n}{n} \sum_{1\leq i < j \leq n} \left[w\left(X_{i,k},X_{i,\ell},Y_i,X_{j,k},X_{j,\ell},Y_j\right) - w_1(X_{i,k},X_{i,\ell},Y_i) - w_1(X_{j,k},X_{j,\ell},Y_j)\right]\, +\,\E\left(U_{n}\right)\\
										&=&\frac{\pi_n}{n} \sum_{1\leq i < j \leq n}\mathcal{H}\left(X_{i,k},X_{i,\ell},Y_i,X_{j,k},X_{j,\ell},Y_j \right),
\end{eqnarray*}
it is therefore clear that
\begin{multline*}
\mathcal{H}\left(X_{i,k},X_{i,\ell},Y_i,X_{j,k},X_{j,\ell},Y_j \right) = w\left(X_{i,k},X_{i,\ell},Y_i,X_{j,k},X_{j,\ell},Y_j\right) - w_1(X_{i,k},X_{i,\ell},Y_i)\\
 - w_1(X_{j,k},X_{j,\ell},Y_j) + \E\left(w\left(X_{i,k},X_{i,\ell},Y_i,X_{j,k},X_{j,\ell},Y_j\right)\right),
\end{multline*}
We see that $\mathcal{H}(\cdot)$ is a symmetric kernel and it's clear that $\E\left[\mathcal{H}\left(X_{i,k},X_{i,\ell},Y_i,X_{j,k},X_{j,\ell},Y_j \right)\right] =0$. 
To calculate $\E\left[\mathcal{H}\left(X_{i,k},X_{i,\ell},Y_i,X_{j,k},X_{j,\ell},Y_j\right)^2\right]$, we will need $\E\left[w\left(X_{i,k},X_{i,\ell},Y_i,X_{j,k},X_{j,\ell},Y_j\right)^2\right]$.
We have
\begin{multline*}
\E\left( w\left(X_{i,k},X_{i,\ell},Y_i,X_{j,k},X_{j,\ell},Y_j\right)^2\right)\\
  \leq 2\E\left[H^2(Y_i) \pi_j^{-2} \gamma_j^2 X_{j,k}^2 X_{j,\ell}^2 \frac{1}{h_j^2}K^2\left(\frac{Y_i-Y_j}{h_j}\right)\right] \\
	+ 2\E\left[H^2(Y_j) \pi_i^{-2} \gamma_i^2 X_{i,k}^2 X_{i,\ell}^2 \frac{1}{h_i^2}K^2\left(\frac{Y_j-Y_i}{h_i}\right)\right].
\end{multline*}
Since 
\begin{eqnarray*}
& &\E\left[\left(H(Y_i) \pi_j^{-1} \gamma_j X_{j,k} X_{j,\ell}\frac{1}{h_j}K\left(\frac{Y_j-Y_i}{h_j}\right)\right)^2\right]\\
& & \hspace{1cm}  = \left(\pi_j^{-1}\gamma_j \right)^2  \frac{1}{h_j} \int \int \left(v_{jk\ell}(uh_j+y_i) - v_{jk\ell}(y_i)\right) H^2(y_i) K^2(u) f(y_i) dy_i du\\
& & \hspace{3cm}  + \left(\pi_j^{-1}\gamma_j \right)^2  \frac{1}{h_j} \int \int v_{jk\ell}(y_i) H^2(y_i) K^2(u) f(y_i) dy_i du\\
& & \hspace{1cm}  \leq  C_1\left(\pi_j^{-1}\gamma_j \right)^2  \int \int  H^2(y_i) |u| K^2(u) f(y_i) dy_i du\\
& & \hspace{3cm}  + C_2\left(\pi_j^{-1}\gamma_j \right)^2  \frac{1}{h_j} \int  v_{jk\ell}(y_i) H^2(y_i)  f(y_i) dy_i \\
& & \hspace{1cm}  = O\left(\frac{\gamma_j^2\pi_j^{-2}}{h_j}\right).
\end{eqnarray*}
It follows that
\begin{equation*}
\E\left( w\left(X_{i,k},X_{i,\ell},Y_i,X_{j,k},X_{j,\ell},Y_j,\right)^2\right)= O\left(\frac{\gamma_i^2\pi_i^{-2}}{h_i} + \frac{\gamma_j^2\pi_j^{-2}}{h_j}\right).
\end{equation*}
From the equality below, we obtain
\begin{eqnarray*}
& &\E\left[\mathcal{H}\left(X_{i,k},X_{i,\ell},Y_i,X_{j,k},X_{j,\ell},Y_j \right)^2\right]\\
& & =  \E\left[ w\left(X_{i,k},X_{i,\ell},Y_{i},X_{j,k},X_{j,\ell},Y_j\right)^2\right] - \left(\E\big[w\left(X_{i,k},X_{i,\ell},Y_{i},X_{j,k},X_{j,\ell},Y_j\right)\big]\right)^2\\
& & = O\left(\frac{\gamma_i^2\pi_i^{-2}}{h_i} + \frac{\gamma_j^2\pi_j^{-2}}{h_j}\right).
\end{eqnarray*}

\medskip

\textbf{Step 4.} : Show that \,$\sqrt{n}\left(U_n - \widehat{U}_n\right) = O_p\left(\sqrt{\frac{\gamma_n}{h_n}}\right)$.
 To achieve this, we only need to compute the convergence rate of $\E\left[\left(U_n - \widehat{U}_n\right)^2\right]$.\\ 
We can already observe that $\E\left(U_n - \widehat{U}_n\right)=0$. Which implies $var\left[\sqrt{n}\left(U_n - \widehat{U}_n\right)\right] = n\E\left[\left(U_n - \widehat{U}_n\right)^2\right]$.\\
We have
\begin{eqnarray*}
\E\left[\left(U_n - \widehat{U}_n\right)^2\right] &=&\left(\frac{\pi_n}{n}\right)^2 \sum_{1\leq i < j \leq n} \E\left(\mathcal{H}\left(X_{i,k},X_{i,\ell},Y_i,X_{j,k},X_{j,\ell},Y_j \right)^2\right)\\
& & + \left(\frac{\pi_n}{n}\right)^2 \sum_{(i_1,j_1, i_2,j_2)\in \sigma_n} \E\bigg(\mathcal{H}\left(X_{i_1,k},X_{i_1,\ell},Y_{i_1},X_{j_1,k},X_{j_1,\ell},Y_{j_1} \right)\\
& & \hspace{3cm} \times\mathcal{H}\left(X_{i_2,k},X_{i_2,\ell},Y_{i_2},X_{j_2,k},X_{j_2,\ell},Y_{j_2} \right) \bigg),
\end{eqnarray*}
where $\sigma_n = \left\{(i_1,j_1, i_2,j_2)\,/\, 1\leq i_1 < j_1 \leq n,\,1\leq i_2 < j_2 \leq n,\, (i_1,j_1)\neq (i_2,j_2)  \right\}$.\\
However, due to the fact that $(i_1,j_1, i_2,j_2) \in \sigma_n$, $\E\big(\mathcal{H}\left(X_{i_1,k},Y_{i_1},X_{j_1,k},Y_{j_1} \right) \times\mathcal{H}\left(X_{i_2,k},Y_{i_2},X_{j_2,k},Y_{j_2} \right) \big) = 0$.\\
In the sequel
\begin{eqnarray*}
& & \left(\frac{\pi_n}{n}\right)^2 \sum_{1\leq i < j \leq n} \E\left(\mathcal{H}\left(X_{i,k},X_{i,\ell},Y_i,X_{j,k},X_{j,\ell},Y_j \right)^2\right) \\
& & \hspace{3cm}\leq  \left(C\frac{\pi_n}{n}\right)^2 \sum_{1\leq i < j \leq n}\left(\frac{\gamma_i^2\pi_i^{-2}}{h_i} + \frac{\gamma_j^2\pi_j^{-2}}{h_j}\right)\\
& & \hspace{3cm}\leq C\left(\frac{\pi_n}{n}\right)^2(n-1)\sum_{i=1}^n\frac{\gamma_i^2\pi_i^{-2}}{h_i}\\
& & \hspace{3cm}\leq C\left(\frac{\gamma_n}{nh_n}\right).
\end{eqnarray*}
Which means 
$$
 \E\left[\left(U_n - \widehat{U}_n\right)^2\right] = O\left(\frac{\gamma_n}{nh_n}\right)\,\mbox{ and thus }\,n\E\left[\left(U_n - \widehat{U}_n\right)^2\right] = O\left(\frac{\gamma_n}{h_n}\right).
$$
Hence the result of this step.

\medskip

\textbf{Step 5.}: Conclusion of the proof\\

From (\ref{eqnc1lem3}), (\ref{eqnlr5a}) and (\ref{eqn1u})
\begin{eqnarray*}
\frac{1}{\sqrt{n}}\sum_{i=1}^n H(Y_i)\widehat{G}_{k\ell,n}(Y_i) &=& \sqrt{n}U_{n} + o_p(1)\\
                                                                 &=& \sqrt{n}\widehat{U}_{n} + o_p(1)\\
																															   &=& \sqrt{n}\sum_{i=1}^n \left[\frac{1}{n}H(Y_i) G_{k\ell}(Y_i) + \pi_n\pi_i^{-1} \gamma_i X_{i,k} H(Y_i) f(Y_i)\right]\\
																																& & \qquad - \,\sqrt{n}\E(H(Y)G_{k\ell}(Y)) + O(\sqrt{n}h_n^4)+ o_p(1).
\end{eqnarray*}
Then
\begin{multline*}
\frac{1}{\sqrt{n}}\sum_{i=1}^n \bigg(H(Y_i)\widehat{G}_{k\ell,n}(Y_i)\, -\, H(Y_i) G_{k\ell}(Y_i)\bigg)\\
 = \sqrt{n}\sum_{i=1}^n\bigg(\pi_n\pi_i^{-1} \gamma_i X_{i,k} X_{i,\ell}  H(Y_i) f(Y_i) \,-\, \frac{1}{n}\E(H(Y)G_{k\ell}(Y))\bigg)\,+\, o_p(1).
\end{multline*}
The proof is completed.\,$\square$

\end{document}